\documentclass[twoside, 11pt, reqno]{amsart}

\usepackage[top=1.1in,bottom=1.1in,left=1.2in,right=1.2in]{geometry}
\usepackage{amsmath,amssymb,amsthm,amsfonts,enumerate, mathrsfs, mathtools, appendix}
\usepackage{relsize}
\usepackage{stmaryrd}
\usepackage{esint}
\usepackage{xcolor}
\usepackage{graphicx}
\usepackage{tikz} 
\usetikzlibrary{patterns,positioning,arrows,shapes,trees}
\usepackage[%colorlinks=false,   
colorlinks=true,
linkcolor=blue,
filecolor=blue,
urlcolor=red,
citecolor=magenta]{hyperref}

\usepackage{fancyhdr}
{} 
\makeatletter
\newcommand{\mylabel}[2]{#2\def\@currentlabel{#2}\label{#1}}
\makeatother

\newcounter{counterConstant} 
\newcommand{\const}[1]{
	\addtocounter{counterConstant}{1}
	\edef#1{\arabic{counterConstant}}
}

\numberwithin{equation}{section}

\numberwithin{equation}{section}

\newtheorem{theorem}{Theorem}[section]
\newtheorem{lemma}[theorem]{Lemma}
\newtheorem{remark}[theorem]{Remark}

\newtheorem{proposition}[theorem]{Proposition}

\def\cD{{\mathcal D}}

\def\cM{{\mathcal M}}

\def\cO{{\mathcal O}}

\def\cS{{\mathcal S}}
\def\cT{{\mathcal T}}

\def\mD{{\mathbb D}}

\def\mH{{\mathbb H}}

\def\mL{{\mathbb L}}

\def\mN{{\mathbb N}}

\def\mQ{{\mathbb Q}}
\def\mR{{\mathbb R}}

\def\R{\mathbb R}

\def\bP{{\mathbf P}}

\def\bT{{\mathbf T}}

\def\bP{{\mathbf P}}
\def\bQ{{\mathbf Q}}
\def\bE{{\mathbf E}}
\def\1{{\mathbf{1}}}

\def\sF{{\mathscr F}}
\def\sG{{\mathscr G}}
\def\sH{{\mathscr H}}

\def\tL{\widetilde{\mathbb{L}}}
\def\tH{\widetilde{\mathbb{H}}}

\def\div{\mathord{{\rm div}}}
\def\eps{\varepsilon}
\def\d{\text{\rm{d}}}

\def\e{\mathrm{e}}

\def\a{\alpha}
\def\om{\omega}
\def\Om{\Omega}

\def\p{\partial}

\def\l{\lambda}
\def\si{\sigma}

\def\<{{\langle}}
\def\>{{\rangle}}

\def\dif{{\mathord{{\rm d}}}}

\def\={&\!\!=\!\!&}

\def\l{\left}
\def\r{\right}

\def\geq{\geqslant}
\def\leq{\leqslant}

\allowdisplaybreaks

\begin{document}
	
\title{SDEs with critical time dependent drifts: strong solutions}
	
\author{Michael R\"ockner and Guohuan Zhao} 
		
\address{Michael R\"ockner: Department of Mathematics, Bielefeld University, Germany \\ and Academy of Mathematics and Systems Science, Chinese Academy of Sciences (CAS), Beijing, 100190, P.R.China}
\email{roeckner@math.uni-bielefeld.de}
	
\address{Guohuan Zhao: Academy of Mathematics and Systems Science, Chinese Academy of Sciences (CAS),
		Beijing, 100190, P.R.China}  
\email{gzhao@amss.ac.cn}

\thanks{Research of Michael and Guohuan is supported by the German Research Foundation (DFG) through the Collaborative Research Centre (CRC) ``Taming uncertainty and profiting from randomness and low regularity in analysis, stochastics and their applications,"-Project-ID 317210226-SFB 1283.}
	
\begin{abstract} 
This paper is a continuation of \cite{rockner2022weak}. Based on a compactness criterion for random fields in Wiener-Sobolev spaces, in this paper, we prove the strong solvability of time-inhomogeneous stochastic differential equations with drift coefficients in critical Lebesgue spaces, which gives an affirmative answer to a longstanding open problem.  As an application, we also prove a regularity criterion for solutions of a stochastic system proposed by Constantin and Iyer (Comm. Pure. Appl. Math. 61(3): 330--345, 2008), which is closely related to the Navier-Stokes equations. 
	\end{abstract}
	
	\maketitle
	\tableofcontents
	
	
	\section{Introduction}
	
	Let $W_t$ be a standard $d$-dimensional Brownian motion on a complete filtered  probability space $(\Om, \sF, \{\sF_t\}_{t\in [0,T]}, \bP)$ and let $b$ be a vector field on $\R^d$ satisfying the following critical Ladyzhenskaya-Prodi-Serrin (LPS) condition: 
	\begin{equation}\label{Con-critical}
	b\in \mL^{p_1}_{q_1}(T):=L^{q_1}([0,T]; L^{p_1}(\R^d)) \mbox{ with } p_1,q_1 \in [2, \infty]  \mbox{ and } \frac{d}{p_1}+\frac{2}{q_1}=1. 
	\end{equation}
	Our primary goal is to solve the following longstanding open problem: does  the  stochastic differential equation (SDE) 
	\begin{equation}\label{Eq-SDE}
	X_{s,t}^x=x+\int_s^t b(r, X_{s,r}^x) \d r+ W_t-W_s, \quad 0\leq s\leq t\leq T, \ x\in \R^d
	\end{equation}
	have a unique strong solution under condition \eqref{Con-critical}? 
	
\subsection{Main result}
Our main result, which gives an almost affirmative answer to the above open problem, reads  as follows. 
\begin{theorem}\label{Th-Main}
    Let $d\geq 3$.  Assume $b$ satisfies one of following two conditions  
    \begin{enumerate}[(a)]
	\item $b\in C([0,T]; L^d)$;   
	\item  $b\in \mL^{p_1}_{q_1}(T)$ with $p_1,q_1\in (2,\infty)$ and $d/p_1+2/q_1=1$.
    \end{enumerate}
    Then \eqref{Eq-SDE} has a strong solution satisfying the following estimate
    \begin{equation}\label{Eq-Krylov}
       \sup_{x\in \R^d}\bE \l(\int_{s}^{T} f(t, X_{s,t}^x) \ \d t \r) \leq C \sup_{z\in \R^d}\|f\1_{B_1(z)}\|_{L^q([s,T]; L^{p})}, 
    \end{equation}
    where $p, q\in (1,\infty)$ with $\frac{d}{p}+\frac{2}{q}<2$ and $C$ is a constant independent with $f$. Moreover, pathwise uniqueness holds within the class of weak solutions satisfying \eqref{Eq-Krylov}, and in this case, the random field $(X_{s,t}^x)_{\substack{x \in \mathbb{R}^d, \\ 0 \leq s \leq t \leq T}}$ satisfies the following properties:
    \begin{enumerate}
	\item 
        for any $r\geq 1$,  
	\begin{equation}\label{Eq-gradient}
		\sup_{\substack{z\in \R^d\\0\leq s\leq t\leq T}} \int_{B_1(z)} \bE |\nabla X_{s,t}^x-\mathrm{I} |^r \d x <\infty. 
	\end{equation}
	\item for all $r\in (d,\infty)$, $\beta\in (0, \frac{1}{2})$, $R>0$, $x_i\in B_R$, $0\leq s_i\leq t_i\leq T$, $i=1,2$  
		\begin{equation}\label{Eq-Holder-X}
		\bE \l|X_{s_1, t_1}^{x_1}-X_{s_2, t_2}^{x_2}\r|^r \leq C \l(|x_1-x_2|^{r-d} + |s_1-s_2|^{\beta (r-d)}+ |t_1-t_2|^{\beta r}\r). 
		\end{equation}
    \end{enumerate}
\end{theorem}
\begin{remark}
    \begin{enumerate}[(i)]
        \item For technical reasons, we need to assume \( d \geq 3 \) and \( p < \infty \) in this paper. The condition \( p < \infty \) can in fact be removed; see \cite{kinzebulatov2025strong}. We also expect that similar conclusions hold for \( d = 2 \), although handling this case appears to be technically more involved.
        \item 
        Our main approach of this paper does not work for the full endpoint case $p_1=d$ and $q_1=\infty$. However, in the later case the  conditional weak well-posedness was proved by the same authors of this article in \cite{rockner2022weak}, provided that the divergence of $b$ satisfies an integrability condition. We conjecture that the strong well-posedness of \eqref{Eq-SDE} holds when $b\in \mL^d_\infty(T)$ and $\div b=0$. 
        \item The main content of this work is the existence, as well as the conditional uniqueness of strong solutions. Whether the restriction \eqref{Eq-Krylov} can be removed for the uniqueness is an interesting and challenging question. 
    \end{enumerate}
\end{remark}

\subsection{Motivation and Previous results}
	
The existence of stochastic flows associated with SDEs with singular drifts and their regularity properties have various applications. For instance, in \cite{flandoli2010well}, using the stochastic characteristics corresponding to \eqref{Eq-SDE}, Flandoli-Gubinelli-Priola studied the existence and uniqueness for the stochastic transport equation in an $L^\infty$-setting, provided that the drift $b$ is $\a$-H\"older continuous uniformly in $t$ and the divergence of $b$ satisfies some integrability condition.  Later, stochastic continuity equations were also considered in \cite{neves2015wellposedness} when $b$ is divergence free and it satisfies the subcritical LPS condition  
	\begin{equation}\label{Con-sub}
	b \in \mL^{p_1}_{q_1}(T)  \mbox{ with } p_1,q_1 \in (2, \infty)  \mbox{ and } \frac{d}{p_1}+\frac{2}{q_1}<1. 
	\end{equation}
	The same SPDEs were also investigated by Fedrizzi-Flandoli in \cite{fedrizzi2013noise},  Mohammed-Nilssen-Proske in  \cite{mohammed2015sobolev} and Beck-Flandoli-Gubinelli-Maurelli in \cite{beck2019stochastic} under different settings (see also the reference therein). 
	
	Our work is also motivated by the deep connection between singular SDEs and Navier-Stokes equations. The velocity field $u$ of an incompressible fluid not subject to an external force in $\R^d$ satisfies the Navier-Stokes equation 
	\begin{subequations}
		\begin{align}
			\p_t u-\tfrac{1}{2}\Delta u+ (\nabla u)u  +\nabla P=&0 \ \mbox{ in } [0,T]\times \R^d,   \label{Eq-NSE1}\\
			\div\, u=&0, \label{Eq-NSE2} \\
			u(0)=&\varphi. \label{Eq-NSE3}
		\end{align}
	\end{subequations}
	The mathematical studies of Navier-Stokes equations have a long history.  In \cite{leray1934mouvement}, Leray considered \eqref{Eq-NSE1}-\eqref{Eq-NSE3} for the initial data $\varphi\in L^2$. He proved that there exists a global in time Leray-Hopf weak solution $u\in\mL^2_\infty(T)$ with $\nabla u\in \mL^2_2(T)$. However, to date, the problem of smoothness of Leray-Hopf weak solutions for the 3D Navier-Stokes equations remains open. Studies by Prodi \cite{prodi1959teorema}, Serrin \cite{serrin1962interior} and Ladyzhenskaya \cite{ladyzhenskaya1967uniqueness} found that the interior smoothness of Leray-Hopf weak solutions to \eqref{Eq-NSE1}-\eqref{Eq-NSE2} is guaranteed, provided that $u\in \mL^{p_1}_{q_1}(T)$, for $p_1\in (d, \infty)$ and $d/p_1+2/q_1\leq 1$ (see also \cite{fabes1972initial} and \cite{giga1986solutions}). These conditional regularity results and their generalizations have culminated with the work of Escauriaza-Seregin-{\v{S}}ver{\'a}k \cite{escauriaza2003backward} and \cite{escauriaza2003solutions} for $p_1=d=3$ and then Dong-Du \cite{dong2009navier} for $p_1=d\geq 3$. On the other hand, in the corresponding Lagrangian description, a fluid particle motion is described by the SDE 
	\begin{subequations}
		\begin{align}
			\d X_t^x &= u(t, X_t^x) \d t +\d W_t, \quad X_0^x=x\in \R^d. \label{Eq-SDE2} 
		\end{align}
		When $u$ is smooth, Constantin-Iyer \cite{constantin2008stochastic} presented an elegant stochastic representation for the solutions to the Navier-Stokes equations, namely 
		\begin{align}
			u(t,x)&=\mathrm{P} \bE \left[\nabla^{\top} (X^x_{t})^{-1} \varphi\left(( X^x_t)^{-1}\right)\right],  \label{Eq-rept}
		\end{align}
	\end{subequations}
	where $\mathrm{P}$ is the Leray projection and $(X_t^x)^{-1}$ is the inverse stochastic flow of \eqref{Eq-SDE2}. Conversely, if $u$ is smooth and $(u,X)$ solves the stochastic system \eqref{Eq-SDE2}-\eqref{Eq-rept}, then $u$ also solves \eqref{Eq-NSE1}-\eqref{Eq-NSE3}. From then on, some researchers started to study \eqref{Eq-NSE1}-\eqref{Eq-NSE3} via investigating the corresponding stochastic Lagrangian paths, see \cite{galeati2025almost}, \cite{rezakhanlou2014regular}, \cite{rezakhanlou2016stochastically}, \cite{zhang2010astochastic} and \cite{zhang2016stochastic} etc. Since the problems of the regularity of solutions to the 3D Navier-Stokes equations are very challenging, two natural questions arise: (1) If the drift term is irregular, does there exist a strong solution \((X^x_t)\) to \eqref{Eq-SDE2} (or \eqref{Eq-SDE}) that admits an invertible, weakly differentiable version, ensuring that the right-hand side of \eqref{Eq-rept} is well-defined? (2) Can one also obtain some conditional regularity results for the stochastic system \eqref{Eq-SDE2}-\eqref{Eq-rept}? 
    
    Our Theorem \ref{Th-Main} establishes that under the assumption that the drift coefficients satisfy the critical LPS conditions, strong solutions to \eqref{Eq-SDE2} indeed exist and has a version that is weakly differentiable with respect to the spatial variable. For the second question above, to simplify our presentation, as in \cite{zhang2010astochastic}, in this paper,  we study the backward stochastic system 
	\begin{equation}\label{Eq-SS}
	\left\{\begin{aligned}
		X_{t, s}^x &=x+\int_{t}^{s} u\left(r,  X_{t, r}^x\right) \mathrm{d} r+\left(W_{s}-W_{t}\right), \quad &-T\leq t\leq s\leq 0 \\
		u(t, x) &=\mathrm{P} \mathbf{E}\left[\nabla^{\top}  X_{t,0}^x \varphi\left(  X_{t, 0}^x \right)\right], \quad &-T\leq t\leq 0 
	\end{aligned}\right.  
	\end{equation}
	corresponding to the backward Navier-Stokes equation 
	\begin{equation}\label{Eq-BNSE}
	\p_t u+\tfrac{1}{2}\Delta u+ (\nabla u)u  +\nabla P=0, \quad \div u=0, \quad u(0)=\varphi
	\end{equation}
	instead of the forward one \eqref{Eq-SDE2}-\eqref{Eq-rept}. With the help of a refined version of  \eqref{Eq-gradient}, we give a regularity criterion for solutions to \eqref{Eq-SS} in Theorem \ref{Th-Main2} below, which can be regarded as an analogue of Serrin's regularity criterion for solutions to the 3D Navier-Stokes equations. 
    
    The study of strong well-posedness of non degenerate It\^o equations with bounded drift coefficients dates back to \cite{zvonkin1974transformation} and  \cite{veretennikov1980strong2}. 
	In \cite{krylov2005strong}, Krylov-R\"ockner obtained the existence and uniqueness of strong solutions to \eqref{Eq-SDE}, when $b$  satisfies the subcritical LPS  condition. After that a number of papers were devoted to generalize the strong well-posedness result, as well as the following gradient estimate for $X$: 
	\begin{equation}\label{Eq-Ugradient}
	\sup_{x\in \R^d} \bE \sup_{t\in[s, T]} |\nabla X_{s,t}^x|^r<\infty, \quad \forall r\geq 1.  
	\end{equation}
	The reader is referred to \cite{fedrizzi2011pathwise}, \cite{lee2017existence}, \cite{rezakhanlou2014regular}, \cite{xia2020lqlp}, \cite{zhang2005strong}, \cite{zhang2011stochastic}, \cite{zhang2016stochastic} and the reference therein for more details. To the best of our knowledge, the strong solvability under the critical condition \eqref{Con-critical} was first touched by Beck-Flandoli-Gubinelli-Maurelli in \cite{beck2019stochastic}, where they proved the pathwise uniqueness to SDE \eqref{Eq-SDE} in a certain class if the initial datum has a diffuse law. Recently, if $b$ belongs to the Orlicz-critical space $L^{q_1,1}([0,T]; L^{p_1})\subsetneqq \mL^{p_1}_{q_1}(T)$ for some $p_1,q_1\in (2,\infty)$ with $d/p_1+2/q_1=1$, by Zvonkin's transformation (cf.  \cite{zvonkin1974transformation}), Nam \cite{nam2020stochastic} showed the existence and uniqueness of strong solutions for SDE \eqref{Eq-SDE}. The key step in using Zvonkin type of change of variables is to construct a homeomorphism by solving the Kolmogorov equation. If, however,  $b$ only satisfies the critical LPS condition \eqref {Con-critical}, this strategy seems impossible to implement. Recently, Krylov \cite{krylov2021strong} proved the strong well-posedness of \eqref{Eq-SDE} for the case that $b(t,x)=b(x)\in L^d(\R^d)$ with $d\geq 3$, which is a significant progress on this topic. His approach is based on an analytical criterion for the existence of strong solutions, originally established in \cite{veretennikov1976explicit}, and also some new estimates obtained in \cite{krylov2021stochastic2} and \cite{krylov2021stochastic1}.  After the authors submitted this manuscript to arXiv, Krylov further advanced the robust methodology originally proposed in \cite{krylov2021strong}. Subsequently, in \cite{krylov2025strong}, he addressed situations where the equation coefficients exhibit time dependence, and in which he discusses in detail conditional and unconditional strong uniqueness. In this study, we adopt a markedly distinct approach from that delineated in \cite{krylov2021strong} and \cite{krylov2025strong}. A succinct overview of this approach will be provided in the following subsection.

    We close this subsection by noting that for \eqref{Eq-SDE}, when \( b \) is a Leray–Hopf solution to the Navier–Stokes equations, Galeati \cite{galeati2025almost} recently proved the uniqueness of both deterministic and stochastic Lagrangian trajectories for Lebesgue almost every initial condition, using an asymmetric Lusin–Lipschitz property of \( u \). However, whether the associated stochastic flow is weakly differentiable remains open. More generally, the weak well-posedness of SDEs with singular drifts far beyond the classical LPS conditions continues to be an active area of research, with notable progress in \cite{grafner2024weak}, \cite{hao2024sdes} and \cite{zhang2020stochastic}, among others.

    \subsection{Approach and Structure}
	The approach in this article is probabilistic, employing ideas from the Malliavin calculus coupled with some estimates for parabolic equations. In \cite{rockner2022weak}, we obtain weak well-posedness of \eqref{Eq-SDE} under a slightly more general condition. So, to get the strong well-posedness, one only needs to show the strong existence due to a fundamental result of Cherny \cite{cherny2002uniqueness}. Our approach to proving strong existence is quite straightforward. Let $\{b_k\}$ be a smooth approximating sequence of the drift $b$ in $\mL^{p_1}_{q_1}(T)$ and $X^x_{s,t}(k)$ be the unique strong solution to \eqref{Eq-SDE} with $b$ replaced by $b_k$. The main effort of the present work is to show that $X^x_{s, t}(k)$ converges to a random field $X^x_{s, t}$, which is a strong solution to equation \eqref{Eq-SDE}. A key ingredient for the convergence of   $X^x_{s, t}(k)$ is the fact that for each $s, t\in [0,T]$ and $R>0$, the sequence $\{X_{s, t}^x(k)\}$ is compact in $L^2(B_R\times \Om)$. The proof for this assertion is based on a compactness criterion for $L^2$ random fields in Wiener spaces (see Lemma \ref{Le-Comp} below or \cite{bally2004relative}) and the following crucial estimate: for any $\a_i\in \{1, 2,\cdots, d\} (i\in \mN_+)$, $n\in \mN_+$ and some $p>1$, 
	\begin{equation}\label{eq-iteration}
	\l\| \bE \int\!\!\!\cdots\!\!\!\int_{s\leq t_1\leq \cdots\leq t_n\leq t}\prod_{i=1}^n \p_{\a_i} f_i(t_i, X_{s, t_i}^x(k)) \ \dif t_1\cdots\dif t_{n}  \r\|_{L^p_{x,loc}} \leq  C^{n+1} \prod_{i=1}^n\|f_i\|_{\mL^{p_1}_{q_1}(s,t)}, 
	\end{equation}
	where $C$ does not depend on $k$ and $\a_i$ (see Lemma \ref{Le-Iteration} below for the precise statement). 
	
	In fact, the framework mentioned above for proving the strong existence of SDEs with {\em bounded} drifts has already been used in  \cite{menoukeu2013variational}, \cite{meyer2010construction} and \cite{mohammed2015sobolev}. However, the main techniques in our paper are essentially  different in comparison with the previous literature. For example, in \cite{mohammed2015sobolev}, using Girsanov's transformation, the main ingredient for the proof of the strong existence result and the gradient estimate for $X$ was reduced to the following estimate: 
	\begin{align*}
		\l\| \bE\int\!\!\!\cdots\!\!\!\int_{0\leq t_1\leq \cdots\leq t_n\leq t} \prod_{i=1}^n  \p_{\a_i} f_i\left(t_{i},x+W_{t_i}\right) \d t_1 \cdots \d t_n \r\|_{L^\infty_x}\leq C^n t^{\frac{n}{2}}(n!)^{-\frac{1}{2}} \prod_{i=1}^n\|f_i\|_{\infty}. 
	\end{align*}
	Such a bound was first obtained by Davie in \cite{davie2007uniqueness} (cf.  \cite[Proposition 2.2]{davie2007uniqueness} and  \cite[Proposition 2.1]{shaposhnikov2016some}) by proving a bound for certain block integrals. Later, Rezakhanlou \cite{rezakhanlou2014regular} also showed that 
	\begin{align*}
		\l\| \bE \int\!\!\!\cdots\!\!\!\int_{0\leq t_1\leq \cdots\leq t_n\leq t} \prod_{i=1}^n  \p_{\a_i} f_i\left(t_{i}, x+W_{t_i}\right) \d t_1 \cdots \d t_n \r\|_{L^\infty_x} \leq C^n t^{\frac{\kappa}{2}} (n!)^{-\frac{\kappa}{2}} \prod_{i=1}^n\|f_i\|_{\mL^p_q(t)}, 
	\end{align*}
	provided that $\kappa:= 1-\frac{d}{p}-\frac{2}{q}>0$. However, when $\kappa=0$, one can not expect to have bounds that are uniform in $x$, and the approach used in \cite{mohammed2015sobolev} and \cite{rezakhanlou2014regular} seems very hard, if not impossible to deal with the critical case ($\kappa=0$). To overcome these essential difficulties, due to the fact that we are in the critical case, in this paper, we  reduce the desired bounded \eqref{eq-iteration} to a uniform in time $L^p$-bound on the solution to a certain parabolic equation with critical drift and a distributional valued  inhomogeneous term (see the discussion before Lemma \ref{Le-Iteration} below). To obtain such uniform bound,  we investigate the PDE mentioned above in Sobolev spaces with mixed norms (see Theorem \ref{Th-Key1} and \ref{Th-Key2}) with the aid of some parabolic versions of Sobolev and Morrey inequalities in mixed norm spaces, which are proved by using Sobolevskii Mixed Derivative Theorem. We note that the PDE results required for proving the main result of this paper can be relaxed, as discussed in \cite{kinzebulatov2025strong}.
	
	\medskip
	
	The rest of this paper is organized as following: In the rest of this section, we list some notations that will be used in this paper frequently. In Section \ref{Sec-Analytic},  we study Kolmogorov equations with  inhomogeneous terms in Sobolev spaces of negative order.  In Section \ref{Sec-Compact}, give a compactness criterion for $L^2$ random fields in Wiener spaces. In Section \ref{Sec-Key}, we derive some crucial uniform estimates for the solutions to certain approximating SDEs. The proof of the main result is presented in Section \ref{Sec-Proof}. In Section \ref{Sec-Application}, we apply our main result to prove a regularity criterion for solutions of a stochastic system, which is closely related to the Navier-Stokes equations.
	
	\subsection{Notations}\label{Subsec-notation}
	We close this section by mentioning some notational conventions used throughout this paper: 
	\begin{itemize} 
		\item $\mN:=\{0, 1,2,\cdots\}$, $\mN_+:= \{1,2,\cdots, \}$. 
		\item The transpose of a matrix $A$ is denoted by $A^{\top}$. 
		\item 
		For a differentiable map $X: \R^d \ni x \mapsto (X^1(x), \cdots, X^{d_1}(x))^{\top}\in \R^{d_1}$, the matrix $\nabla X(x)$ is defined by
		$$
		\nabla X(x)= \l( \begin{matrix}
			\p_1X^1(x)& \p_2 X^1(x) & \cdots & \p_d X^1(x)\\
			\p_1X^2(x)& \p_2 X^2(x) & \cdots & \p_d X^2(x)\\ 
			\cdots& \cdots & \cdots & \cdots\\
			\p_1X^{d_1}(x)& \p_2 X^{d_1}(x) & \cdots & \p_d X^{d_1}(x)
		\end{matrix} \r). 
		$$
		\item Given $S, T\in [-\infty,\infty]$, set 
		$$
		\Delta_n(S,T):= \{(t_1,\cdots, t_n)\in \R^n: S\leq t_1\leq \cdots\leq t_n\leq T\}, \quad \Delta_n(T):= \Delta_n(0,T). 
		$$
		\item Assume that for each $i\in\{1,2\}$, $(X_i, \Sigma_i, \mu_i)$ is a measure space. Suppose that $f: X_1\times X_2\to \R$, define 
		$$
		\|f\|_{L^{p_1}_{x_1}(\mu_1)L^{p_2}_{x_2}(\mu_2)}:= \l[\int_{X_1} \l( \int_{X_2} |f(x_1, x_2)|^{p_2} \mu_2(\d x_2) \r)^{1/p_2} \mu_1(\d x_1) \r]^{1/p_1}. 
		$$
		\item For each $p,q\in[1,\infty]$, the space $L^q([S,T];L^p(\mR^d))$ is denoted by $
		\mL^p_q(S, T)$. For any $p,q\in (1,\infty), s\in \mR$, define $\mH^{s, p}_q(S, T)=L^q([S,T]; H^{s, p}(\mR^d))$, 
		where $H^{s, p}=(1-\Delta)^{-s/2}L^p$ is the Bessel potential space.  
		
		\item 
		Throughout this paper, we fix a cutoff function 
		\begin{align*}
			\chi\in C^\infty_c(\mR^{d};[0,1]) \mbox{ with $\chi|_{B_1}=1$ and $\chi|_{B^c_2}=0$}. 
		\end{align*}
		For $r>0$ and $x\in\mR^{d}$, let $\chi^{z}_r(x):=\chi\l(\frac{x-z}{r}\r)$. For any $p,q\in [1,\infty]$, define 
		$$
		\widetilde L^p:= \l\{ f\in L^p_{loc}(\R^d): \|f\|_{\widetilde L^p} := \sup_{z\in \R^d} \|f\chi_1^z\|_{p} <\infty\r\}
		$$
		and 
		\begin{align}\label{Def-loc-Lpq}
			\widetilde \mL^p_q(S,T):= \l\{f\in L^q([S,T]; L^p_{loc}(\mR^d)): \|f\|_{\widetilde \mL^p_q(S, T)}:= \sup_{z\in \mR^d} \| f\chi^z_1 \|_{\mL^p_q(S, T)}<\infty.  \r\}. 
		\end{align}
		The localized Bessel potential space is defined as follows: 
		\begin{align*}
			\widetilde \mH^{s,p}_q(S, T):= \left\{f\in L^q([S,T]; H^{s, p}_{loc}(\R^d)): \|f\|_{\widetilde \mH^{s,p}_q(S, T)}:= \sup_{z\in \R^d}\|f\chi^z_1\|_{\mH^{s,p}_q(S, T)} <\infty\right\}. 
		\end{align*}
		\item For simplicity, we set 
		$$\mL^{p}_q(T):= \mL^{p}_q(0,T), \ \mL^{p}_q:= L^q(\R; L^p),\  \mH^{s, p}_q(T):= \mH^{s, p}_q(0,T),\  \mH^{s, p}_q:= L^q(\R; H^{s,p}) 
		$$ 
		and 
		$$
		\tL^{p}_q(T)=\tL^p_q(0,T), \quad \tH^{s, p}_q(T):= \tH^{s, p}_q(0,T).
		$$
	\end{itemize}

	\section{Some auxiliary analytic results}\label{Sec-Analytic}
	In this section, we study the Kolmogorov equations with  inhomogeneous terms in localized Sobolev spaces of negative order. By the basic localization procedure, it is sufficient to study the equations in the usual Sobolev spaces (see Remark \ref{Rek-local-W2}). These  analytic results, which are of their own interest, will play a crucial role in proofs for the main results. 
	
	The following conclusions are variants of Theorem 1.1 and 1.2 in  \cite{krylov2001heat}.  
	\const{\CH2}
	\begin{lemma}\label{Le-heat}
	Let $p,q\in (1,\infty)$ and $\a \in \R$. 
	\begin{enumerate}
		\item Assume $\lambda>0$, $\mu\geq 0$. For each $u\in L^q(\R; H^{\a+2,p})\cap H^{1,q}(\R; H^{\a,p})$, 
		\begin{equation}\label{Eq-CZest}
		\|\p_t u\|_{\mH^{\a, p}_q} + \lambda \|\nabla^2 u\|_{\mH^{\a, p}_q} + \mu\|u\|_{\mH^{\a, p}_q} \leq C \|(\p_t-\lambda \Delta+\mu) u\|_{\mH^{\a, p}_q}, 
		\end{equation}
		where $C$ only depends on $d, p, q$. 
		
		\item Assume that $f\in  \mH^{\a, p}_q(T)$,  then the following heat equation admits a unique solution in $\mH^{\a+2,p}_q(T)$: 
		$$
		\p_t u-\tfrac{1}{2}\Delta u= f \ \mbox{ in } (0,T)\times \R^d, \quad u(0)=0. 
		$$
		Moreover, 
		\begin{equation}\label{Eq-Heat-est}
		\|\p_t u\|_{\mH^{\a, p}_q(T)}+\|u\|_{\mH^{\a+2, p}_q(T)}\leq C_{\CH2} \|f\|_{\mH^{\a, p}_q(T)}, 
		\end{equation}
		where $C_{\CH2}$ only depends on $d, p, q,T$. 
	\end{enumerate}
	\end{lemma}
	Consider the following Kolmogorov equation associated with \eqref{Eq-SDE}: 
	\begin{equation}\label{Eq-Kolmogorov}
	\p_t u = \tfrac{1}{2}\Delta u + b\cdot \nabla u+f, \quad u(0)=0. 
	\end{equation}
	
	Throughout this paper, we fix a smooth function $\rho\in C_c^\infty(\R^d)$ satisfying $\rho\geq 0$ and $\int \rho=1$, and set $\rho_m(\cdot):= m^d\rho(m\cdot)$. 
	\subsection{Case (a): \texorpdfstring{$b\in C([0,T]; L^d(\R^d))$}{}}  For any $f\in \mL^d_\infty(T)$, define  
	\begin{equation}\label{Eq-Kf}
	K_f(m):= \sup_{t\in [0,T]} \|f(t)-f(t)*_x\rho_m\|_{L^d}. 
	\end{equation}
	\begin{proposition}\label{Prop-CLd}
	Suppose that $f\in C([0,T]; L^d)$, then  $K_f(m)\to 0$ as $m\to\infty$. 
	\end{proposition}
	\begin{proof}
	Since the map $f:[0,T]\to L^d$ is uniformly continuous, for each $\eps>0$ there is a constant $\delta>0$ such that 
	$$
	\sup_{\substack{t_1,t_2\in [0,T]; \\|t_1-t_2|\leq \delta}}\|f(t_1)-f(t_2)\|_{L^d}<\eps/2.
	$$ 
	Assume that $k=\{0, 1, 2,\cdots, [T/\delta]\}$ and $t\in [k\delta, (k+1)\delta\wedge T]$, then 
	\begin{align*}
		&\limsup_{m\to\infty}\|f(t)-f_m(t) \|_{L^d} \\
		\leq& \|f(t)-f(k\delta)\|_{L^d} + \limsup_{m\to\infty}\|f(k\delta)-f_m(k\delta)\|_{L^d}+ \limsup_{m\to\infty}\|[f(k\delta)-f(t)]*_x\rho_m\|_{L^d}\\
		\leq& 2 \| f(t)-f(k\delta) \|_{L^d}\leq 2\sup_{\substack{t_1,t_2\in [0,T]; \\|t_1-t_2|\leq \delta}}\|f(t_1)-f(t_2)\|_{L^d}<\eps. 
	\end{align*}
	Thus, $\lim_{m\to\infty}K_f(m)=0$. 
	\end{proof}
	The following theorem will plays a crucial role in the proof of the first case of our main result. 
	\const{\CW2}
	\begin{theorem}\label{Th-Key1}
	Let $d\geq 3$, $\a\in \{0,-1\}$ 
	and  $\{a(m)\}_{m\in \mN_+}$ be a sequence converging to zero.  Assume $b\in \mL^d_\infty(T)$ and $K_b(m)\leq a(m)$. Suppose that 
	$p\in (1,d)$ and $q\in (1,\infty)$ if $\a=0$, or $p\in (d/(d-1),d)$ and $q\in (1,\infty)$ if $\a=-1$. Then for any $f\in \mH^{\a,p}_q(T)$, equation \eqref{Eq-Kolmogorov} admits a solution in $\mH^{\a+2,p}_q(T)$. Moreover, 
	\begin{equation}\label{Eq-W2}
	\sup_{t\in (0,T]}t^{-1} \|u\|_{\mH^{\a,p}_q(t)}+\|\p_t u\|_{\mH^{\a,p}_q(T)} + \|u\|_{\mH^{\a+2,p}_q(T)} \leq C_{\CW2}\|f\|_{\mH^{\a,p}_q(T)}, 
	\end{equation}
	where $C_{\CW2}$ only depends on $d, p,q, T, \|b\|_{\mL^d_\infty(T)}, \{a(m)\}$. 
	\end{theorem}
	\begin{proof}
	Below we only give the proof for the case that $\a=-1$ (the case $\a=0$ is simpler). To prove the desired result, it suffices to show \eqref{Eq-W2} assuming that the solution already exists, since the method of continuity is applicable. 
	Let $b_m= b*_x\rho_m$ and $\bar b_m= b-b_m$. Noting that $dp/(p+d)>1$, by Sobolev embedding and H\"older's inequality, we have 
	\const{\CSobolev1}
	\begin{equation}\label{eq-barbmu'}
	\begin{aligned}
		\|\bar b_m\cdot \nabla u\|_{\mH^{-1, p}_q(t)} \leq& C_{\CSobolev1} \|\bar b_m\cdot \nabla u\|_{\mL^{\frac{dp}{p+d}}_q(t)} \leq C_{\CSobolev1} \|\bar b_m\|_{\mL^{d}_\infty(t)} \|\nabla u\|_{\mL^p_q(t)}\\
		\leq & C_{\CSobolev1} a(m) \| u\|_{\mH^{1,p}_q(t)},  
	\end{aligned}
	\end{equation}
	where $t\in [0,T]$ and $C_{\CSobolev1}$ only depends on $d, p$. Similarly, 
	\begin{equation}\label{eq-bmu'}
	\begin{aligned}
		\|b_m \cdot \nabla u\|_{\mH^{-1, p}_q(t)} \leq & \|\div (b_m \, u)\|_{\mH^{-1, p}_q(t)}+\|\div b_m \cdot u\|_{\mL^{p}_q(t)} \\
		\leq &C \|b_m\|_{L^\infty([0,T]; C_b^1)} \|u\|_{\mL^{p}_q(t)}\\
		\leq& C \l( \|\rho_m \|_{d/d-1}+  \|\nabla \rho_m \|_{d/d-1}\r) \|b\|_{\mL^d_\infty(T)}\|u\|_{\mL^{p}_q(t)}\\
		\leq& C m^2 \|u\|_{\mL^{p}_q(t)},  
	\end{aligned}
	\end{equation}
	where $t\in [0,T]$ and $C$ only depends on $d, p, \|b\|_{\mL^d_\infty(T)}$. Thanks to Lemma \ref{Le-heat}, for each $t\in [0,T]$, 
	\begin{align*}
		&\|\p_t u\|_{\mH^{-1,p}_q(t)} + \|u\|_{\mH^{1,p}_q(t)}\\
		\leq &  C_{\CH2} \l(\|b\cdot\nabla u\|_{\mH^{-1,p}_q(t)} +\|f\|_{\mH^{-1,p}_q(t)}\r)\\
		\overset{\eqref{eq-barbmu'}, \eqref{eq-bmu'}}{\leq}& C_{\CH2}\l( C_{\CSobolev1} a(m) \|\nabla u\|_{\mL^p_q(t)} +C m^2 \| u\|_{\mL^p_q(t)}+ \|f\|_{\mH^{-1,p}_q(t)}\r). 
	\end{align*}
	Letting $m$ be large enough such that $C_{\CH2}C_{\CSobolev1} a(m)\leq 1/2$ and using interpolation, we obtain 
	\begin{equation}\label{eq-W2}
	I(t):=\|\p_t u\|_{\mH^{-1,p}_q(t)}^q + \|u\|_{\mH^{1,p}_q(t)}^q \leq C \l( \|u\|_{\mH^{-1, p}_q(t)}^q + \|f\|_{\mH^{-1,p}_q(t)}^q\r),
	\end{equation}
	where $t\in [0,T]$ and $C$ only depends on $d, p, q, T, \|b\|_{\mL^d_\infty(T)}, \{a(m)\}$. 
	One the other hand, recalling that $u(0)=0$ and using H\"older's inequality, we have 
	\begin{equation}\label{eq-H0-1}
	\begin{aligned}
		\| u\|_{\mH^{-1, p}_{q}(t)}^{q}=&\int_0^t \|u(\tau,\cdot)\|_{H^{-1, p}}^{q} \d \tau=  
		\int_0^t \l\|\int_0^\tau\p_t u(\sigma,\cdot) \d \sigma\r\|_{H^{-1, p}}^{q} \d \tau\\
		\leq & \int_0^t  \tau^{q-1} \|\p_t u\|_{ \mH^{-1, p}_{q}(\tau)}^{q}  \d \tau \leq \l[ T^{q-1}\int_0^t I(\tau) \d \tau\r]  \wedge \l[ \frac{t^q}{q}  \|\p_t u\|_{ \mH^{-1, p}_{q}(t)}^{q} \r]. 
	\end{aligned}
	\end{equation}
	Combining this and \eqref{eq-W2}, we obtain 
	$$
	I(t) \leq C \|f\|_{\mH^{-1,p}_q(T)}^q + C \int_0^t I(\tau) \d \tau. 
	$$
	Grönwall's inequality yields, 
	\begin{equation}\label{eq-H2}
	\|\p_t u\|_{\mH^{-1,p}_q(T)} + \|u\|_{\mH^{1,p}_q(T)} \leq  C \|f\|_{\mH^{-1,p}_q(T)}.  
	\end{equation}
	Noting that \eqref{eq-H0-1} also implies 
	$$
	\sup_{t\in (0,T]} t^{-1}\|u\|_{\mH^{-1,p}_q(t)} \leq C(q) \|\p_t u\|_{\mH^{-1,p}_q(T)}, 
	$$
	together with \eqref{eq-H2}, we obtain \eqref{Eq-W2}. 
	$$
	\sup_{t\in (0,T]} t^{-1}\| u\|_{\mH^{-1, p}_q(t)}+\|\p_t u\|_{\mH^{-1,p}_q(T)} + \|u\|_{\mH^{1,p}_q(T)} \leq  C_{\CW2} \|f\|_{\mH^{-1,p}_q(T)},   
	$$
	where $C_{\CW2}$ only depends on $d, p,q ,T, \|b\|_{\mL^d_\infty(T)}, \{a(m)\}$. So,  we complete our proof. 
	\end{proof}
	
	\subsection{Case (b): \texorpdfstring{$b\in \mL^{p_1}_{q_1}(T)$}{} with \texorpdfstring{$p_1, q_1\in (2,\infty)$}{} and \texorpdfstring{$d/p_1+2/q_1=1$}{}}  In this case, to obtain a result similar to Theorem \ref{Th-Key1}, we need to prove some parabolic Morrey and Sobolev inequalities. This can be achieved by using the Mixed Derivative Theorem, which goes back to the work of Sobolevskii  (cf. \cite{sobolevskii1977fractional}). 
	
	Let $X$ be a Banach space and let $A: D(A) \to X$ be a closed, densely defined linear operator with dense range. Then $A$ is called sectorial, if 
	$$
	(0, \infty) \subseteq \rho(-A) \quad \text { and } \quad\left\|\lambda(\lambda+A)^{-1}\right\|_{X\to X} \leq C, \quad \lambda>0,  
	$$
	where $\rho(-A)$ is the resolvent set of $-A$. 
	Set 
	$$
	\Sigma_{\phi}:=\{z \in \mathbb{C} \backslash\{0\}:|\arg z|<\phi\}. 
	$$
	We recall that 
	$$
	\phi_{A}:=\inf \left\{\phi \in[0, \pi): \Sigma_{\pi-\phi} \subseteq \rho(-A), \sup _{z \in \Sigma_{\pi-\phi}}\left\|z(z+A)^{-1}\right\|_{X\to X}<\infty\right\}
	$$
	is the the spectral angle of $A$. For each $\theta\in (0,1)$, define 
	$$
	A^\theta x : = \frac{\sin \theta \pi}{\pi} \int_0^\infty \lambda^{\theta-1}  (\lambda + A)^{-1} Ax \ \d \lambda,  \quad x\in D(A) 
	$$
	and 
	$$
	A^{-\theta} x:=  \frac{\sin \theta \pi}{\pi} \int_0^\infty \lambda^{-\theta}  (\lambda + A)^{-1} x \ \d \lambda,   \quad x\in X. 
	$$
	
	We need the following Sobolevskii Mixed Derivative Theorem (cf. \cite{sobolevskii1977fractional}). 
	\begin{lemma}[Mixed Derivative Theorem]\label{Le-MixD}
	Let $A$ and $B$ be two sectorial operators  in a Banach space $X$ with spectral angles $\phi_A$ and $\phi_B$, which are commutative  and satisfy the parabolicity condition $\phi_A+\phi_B<\pi$. Then the coercivity estimate 
	$$
	\|A x\|_{X}+\lambda \|B x\|_{X} \leq M\|A x+\lambda B x\|_{X}, \quad \forall x \in D(A) \cap D(B), \, \lambda >0
	$$
	implies that 
	$$
	\left\|A^{(1-\theta)} B^{\theta} x\right\|_{X} \leq C\|A x+B x\|_{X}, \quad \forall x \in D(A) \cap D(B), \, \theta \in[0,1], 
	$$
	\end{lemma}
	The following parabolic type Sobolev and Morrey inequalities will be used frequently in this work. 
	\const{\CPSobolevone}
	\const{\CPSobolevtwo}
	\const{\CPMorrey}
	\begin{lemma}\label{Le-Inter}
	Let $p, q\in (1,\infty)$, $r\in(p,\infty)$, $s\in (q,\infty)$ and $\a\in \R$. Assume $\p_t u\in \mH^{\a, p}_q(T)$, $u\in \mH^{\a+2,p}_q(T)$ and $u(0)=0$.  
	\begin{enumerate}
		\item  If $1<d/p+2/q=d/r+2/s+1$, then 
		\begin{equation}\label{Eq-PSobolev1}
		\|u\|_{\mH^{\a+1, r}_s(T)}\leq C_{\CPSobolevone} \l( \|\p_{t} u\|_{\mH_{q}^{\a, p}(T)} + \|u\|_{\mH_{q}^{\a+2, p}(T)}\r), 
		\end{equation}
		where $C_{\CPSobolevone}$ only depends on $d, p, q,r,s$.
		\item If $2<d/p+2/q=d/r+2/s+2$, then
		\begin{equation}\label{Eq-PSobolev2}
		\|u\|_{\mH^{\a, r}_s(T)} \leq C_{\CPSobolevtwo} \l( \|\p_{t} u\|_{\mH_{q}^{\a, p}(T)} + \|u\|_{\mH_{q}^{\a+2, p}(T)}\r), 
		\end{equation}
		where $C_{\CPSobolevtwo}$ only depends on $d, p, q, r,s$.  
		\item If $0\leq \theta<1-1/q$, for any $t_1,t_2\in [0,T]$, 
		\begin{equation}\label{Eq-PMorrey}
			\|u(t_1)-u(t_2)\|_{H^{\a+2\theta, p}} \leq C_{\CPMorrey} |t_1-t_2|^{1-1/q-\theta} \l(\left\|\p_{t} u\right\|_{\mH_{q}^{\a, p}(T)} + \|u\|_{\mH_{q}^{\a+2, p}(T)}\r), 
		\end{equation}
		where $C_{\CPMorrey}$ only depends on $d, p, q, \theta$.
	\end{enumerate}
	\end{lemma}
	\begin{proof}
	By considering $(1-\Delta)^{\a/2} u$ instead of $u$, we see that without loss of generality
	we may assume $\a=0$.  
	Let $X= L^q(\R; L^p(\R^d))$, $A=1+\p_t$ and $B=1-\Delta$ in Lemma \ref{Le-MixD}. It is well-known that 
	$$
	\phi_A=\frac{\pi}{2} \ \mbox{ and }\  \phi_B=0. 
	$$
	Due to \eqref{Eq-CZest}, for all $\lambda>0$ we have  
	\begin{align*}
		\|A u\|_X+\lambda \|B u\|_X=&
		\|u+\p_t u\|_{\mL^p_q}+\lambda \| u- \Delta u\|_{\mL^p_q}\\
		\leq& C \l( \|\p_tu\|_{\mL^p_q}+ \lambda \|\nabla^2 u\|_{\mL^p_q}+ (1+\lambda) \|u\|_{\mL^p_q} \r) \\
		\leq & C \|(u+\p_t u)+ \lambda(u-\Delta u)\|_{\mL^p_q}= C \|A u+\lambda B u\|_X,  
	\end{align*}
	where $C$ only depends on $d, p,q$. 
	Thanks to Lemma \ref{Le-MixD}, we obtain 
	\begin{align}\label{eq-AB}
		\|A^{1-\theta}B^{\theta} u\|_{\mL^p_q} \leq C \|\p_t u-\Delta u+2u \|_{\mL^p_q}\leq C\l( \|\p_t u\|_{\mL^p_q}+ \|u\|_{\mH^{2,p}_q} \r), 
	\end{align}
	for all $u\in H^{1,q}\left(\mathbb{R}, L^{p}(\R^d)\right) \cap L^{q}\left(\mathbb{R}, H^{2,p}(\R^d)\right)$.  For any $\a\in (0,1), q\in (1,\infty)$ and $f\in L^q(\R)$, we have 
	$$
	\sF(  (1-\p_{tt}^2)^{\a/2}(1+\p_t)^{-\a} f) = \frac{(1+4\pi^2 |\xi|^2)^{\a/2}}{(1+\mathrm{i} 2\pi \xi)^\a} \sF(f)(\xi) =:m(\xi) \sF(f)(\xi), 
	$$
	where $\sF(f)(\xi):=\int_{u} \e^{\mathrm{i}2\pi x\xi} f(u)\d u$ is the Fourier transform of $f$. 
	Since $|\xi|^k m^{(k)}(\xi)\leq C_k<\infty$,  by Mikhlin's multiplier theorem, the operator $ (1-\p_{tt}^2)^{\a/2}(1+\p_t)^{-\a}$ is bounded on $L^q(\R)$. Therefore, 
	\begin{equation}\label{eq-pinter}
	\begin{aligned}
		\|u\|_{H^{1-\theta, q}(\R; H^{2\theta, p}(\R^d))}=&\|(1-\p_{tt}^2)^{\frac{1-\theta}{2}}(1-\Delta)^{\theta} u\|_{\mL^p_q} \leq C \|A^{1-\theta} B^\theta u \|_{\mL^p_q}\\\overset{\eqref{eq-AB}}{\leq}& C\l( \|\p_t u\|_{\mL^p_q}+ \|u\|_{\mH^{2,p}_q} \r), \quad \forall \theta\in [0,1].
	\end{aligned} 
	\end{equation}
	If $u\in \mH^{2,p}_q(T)$, $\p_t u\in\mL^p_q(T)$ and $u(0,x)=0$, we extend $u$ by 
	$$
	\bar{u}(t,x):=\left\{\begin{array}{ll}
		u(t, x) & \text { if } t\in [0, T]\\
		-3 u\left(2T-t, x\right)+4 u\left(\frac{3T}{2}-\frac{t}{2}, x\right) & \text { if } t\in [T, 2T]\\
		4 u\left(\frac{3T}{2}-\frac{t}{2}, x\right) & \text { if } t\in [2T, 3T]\\
		0& \text { othewise. } 
	\end{array}\right.
	$$
	By the definition of $\bar u$, one sees that 
	\begin{equation}\label{eq-baru-u}
	\|\p_t \bar u\|_{\mL^p_q}+\|\bar u\|_{\mH^{2,p}_q} \leq C\l(  \|\p_t u\|_{\mL^p_q(T)}+\|u\|_{\mH^{2,p}_q(T)} \r). 
	\end{equation}
	Letting $\theta=\frac{1}{2}+\frac{d}{2p}-\frac{d}{2r}=1+\frac{1}{s}-\frac{1}{q}\in [\frac{1}{2},1]$, the 
	Sobolev inequality and the above estimates imply
	\begin{align*}
		\|u\|_{\mH^{1, r}_s(T)}\leq& \|\bar u\|_{\mH^{1,r}_s} \leq \|\bar u\|_{H^{1-\theta, q}(\R; H^{2\theta, p})}\\
		\overset{\eqref{eq-pinter}}{\leq}& C \l( \|\p_{t} \bar u\|_{\mL_{q}^{p}}+ \|\bar u\|_{\mH_{q}^{2, p}}\r) \overset{\eqref{eq-baru-u}}{\leq}  C_{\CPSobolevone} \l(  \|\p_t u\|_{\mL^p_q(T)}+\|u\|_{\mH^{2,p}_q(T)} \r). 
	\end{align*}
	So, we complete our proof for \eqref{Eq-PSobolev1}. \eqref{Eq-PSobolev2} can be proved similarly. 
	
	For \eqref{Eq-PMorrey}, if $\theta<1-1/q$, by Morrey's inequality, we have 
	\begin{align*}
		\sup_{t_1, t_2\in [0,T]}\frac{\|u(t_1)-u(t_2)\|_{H^{2\theta, p}}}{|t_1-t_2|^{1-1 / q-\theta}} \leq& C \|\bar u\|_{H^{1-\theta, q}(\R; H^{2\theta, p})} \\
		\leq& C_{\CPMorrey} \l(\left\|\p_{t} u\right\|_{\mL_{q}^{p}(T)}+\|u\|_{\mH_{q}^{2, p}(T)} \r). 
	\end{align*}
	So, we complete our proof. 
	\end{proof}
	\begin{remark}\label{Rek-Inter}
	It is not hard to see that if $\p_t u\in \tH^{\a, p}_q(T)$, $u\in \tH^{\a+2,p}_q(T)$ and $u(0)=0$, then all the conclusions in Lemma \ref{Le-Inter} still hold if $\mH^{\cdots}_{\cdots}$ and $H^{\cdots}$ are replaced by $\tH^{\cdots}_{\cdots}$ and $\widetilde{H}^{\cdots}$, respectively.
	\end{remark}
	
	For any $f\in \mL^{p_1}_{q_1}(T)$, set 
	\begin{equation}\label{Eq-K'f}
	K'_f(m):=\| f-f\1_{\{|f|\leq m\}}\|_{\mL^{p_1}_{q_1}(T)}
	\end{equation}
	and 
	\begin{equation}\label{eq-om-f}
	\om_f(\delta):= \sup_{0\leq S\leq T-\delta} \|f\|_{\mL^{p_1}_{q_1}(S, S+\delta)}. 
	\end{equation}
	Since $q_1<\infty$, one sees that $K'_f(m) \to 0, \ \mbox{ as } m\to \infty$ and $\om_f(\delta)\to 0$, as $\delta\to 0$. 
	
	\medskip
	
	Next we give an analogue of Theorem \ref{Th-Key1}, which is crucial in the proof of the second case of Theorem \ref{Th-Main}.
	\begin{theorem}\label{Th-Key2}
	Let $d\geq 3$, $p_1,q_1\in (2,\infty)$ with $d/p_1+2/q_1=1$, and $\{a(m)\}_{m\in \mN_+}$ be a sequence converging to zero. Assume that $b\in \mL^{p_1}_{q_1}(T)$ and $K'_b(m)\leq a(m)$, 
	\begin{enumerate}
		\item if $p\in (1,p_1)$ and $q\in (1,q_1)$,  then for any $f\in \mL^{p}_q(T)$, equation \eqref{Eq-Kolmogorov} admits a solution $u$ in $\mH^{2,p}_q(T)$ and 
		\begin{equation}\label{Eq-W2-u}
		\|\p_t u\|_{\mL^{p}_q(T)} + \|u\|_{\mH^{2,p}_q(T)} \leq C\|f\|_{\mL^p_q(T)},
		\end{equation}
		where $C$ only depends on $d,p_1, q_1, p,q, T, \{a(m)\}$ and is increasing in $T$;  
		\item if $p\in (p_1/(p_1-1),p_1)$ and $q\in (q_1/(q_1-1),q_1)$, then for any $f\in \mH^{-1,p}_q(T)$, equation \eqref{Eq-Kolmogorov} admits a solution $u$ in $\mH^{1,p}_q(T)$, and $u=v+w$ with $v, w$ satisfying 
		\begin{equation}\label{Eq-W2-v}
		\|\p_t v\|_{\mH^{-1,p}_q(T)} + \|v\|_{\mH^{1,p}_q(T)} \leq C \|f\|_{\mH^{-1,p}_q(T)},
		\end{equation}
		and 
		\begin{equation}\label{Eq-W2-w}
		\|\p_t w\|_{\mL^{p'}_{q'}(T)} + \|w\|_{\mH^{2,p'}_{q'}(T)} \leq C\|f\|_{\mH^{-1,p}_q(T)},  
		\end{equation}
		where $p'=\frac{p_1p}{p_1+p}>1$,  $q'=\frac{q_1q}{q_1+q}>1$,  and $C$ only  depends on $d, p_1, q_1, p, q, T, \{a(m)\}$ and is increasing in $T$. 
	\end{enumerate}
	\end{theorem}
	\begin{proof}
	To prove the desired result, we only need to prove \eqref{Eq-W2-u}, \eqref{Eq-W2-v} and \eqref{Eq-W2-w} assuming that the solution already exists, since the method of continuity is applicable.  
	
	(1). Let $b_m:= b\1_{\{|b|\leq m\}}$. Rewrite \eqref{Eq-Kolmogorov} as 
	$$
	\p_t u-\tfrac{1}{2}\Delta u=f+ b_m\cdot\nabla u+(b-b_m)\cdot \nabla u. 
	$$
	Thanks to Lemma \ref{Le-heat}, for any $t\in [0,T]$ we have 
	\begin{align*}
		\begin{aligned}
			\|\p_t u\|&_{ \mL^{p}_{q}(t)} + \|u\|_{ \mH^{2, p}_{q}(t)} \\
			&\leq C_{\CH2}  \l( \|f\|_{\mL^{p}_{q}(t)}+ m \|\nabla u\|_{ \mL^{p}_{q}(t)}+
			\|(b-b_m)\cdot \nabla u\|_{ \mL^{p}_{q}(t)} \r), 
		\end{aligned}
	\end{align*}
	where $C_{\CH2}=C_{\CH2}(d, p, q, T)$. Letting $1/r=1/p-1/p_1$ and $1/s=1/q-1/q_1$, by \eqref{Eq-PSobolev1} we have 
	\begin{align*}
		\begin{aligned}
			\|(b-b_m)\cdot \nabla u\|_{ \mL^{p}_{q}(t)} \leq& \|(b-b_m)\|_{\mL^{p_1}_{q_1}(t)} \|\nabla u\|_{\mL^r_s(t)}\\
			\overset{\eqref{Eq-PSobolev1}}{\leq}& C_{\CPSobolevone}a(m) \l(\|\p_t u\|_{\mL^p_q(t)}+\|u\|_{\mH^{2, p}_{q}(t)}\r). 
		\end{aligned}
	\end{align*}
	We choose $m=N$, which is sufficiently large so that $C_{\CH2}C_{\CPSobolevone}a(N)\leq 1/2$. Therefore, 
	\begin{align}\label{eq-It}
		\begin{aligned}
			I(t):=& \|\p_t u\|_{ \mL^{p}_{q}(t)}^{q} +\|u\|_{\mH^{2,p}_{q}(t)}^{q} \leq C\l( \|f\|_{ \mL^{p}_{q}(t)}^{q}+ N^{q} \|\nabla u\|_{ \mL^{p}_{q}(t)}^{q}\r). 
		\end{aligned}
	\end{align}
	Noting that 
	\begin{equation}\label{eq-u-I}
	\begin{aligned}
		\| u\|_{\mL^{p}_{q}(t)}^{q}=&\int_0^t \|u(\tau,\cdot)\|_{L^{p}}^{q} \d \tau=  
		\int_0^t \l\|\int_0^\tau\p_t u(\sigma,\cdot) \d \sigma\r\|_{L^{p}}^{q} \d \tau\\
		\leq & \int_0^t  \tau^{q-1} \|\p_t u\|_{ \mL^{p}_{q}(\tau)}^{q}  \d \tau \leq C(T, q)\int_0^t I(\tau) \d \tau, 
	\end{aligned}
	\end{equation}
	and using an interpolation inequality, we obtain 
	\begin{align}\label{eq-inter}
		\begin{aligned}
			\|\nabla u\|_{ \mL^{p}_{q}(t)}^{q} \leq &\delta \|\nabla^2 u\|_{ \mL^{p}_{q}(t)}^{q} + C_\delta\| u \|_{ \mL^{p}_{q}(t)}^{q}\\
			\leq &\delta I(t) + C_\delta \int_0^t I(\tau) \d \tau, \quad (\forall \eps>0). 
		\end{aligned}
	\end{align}
	Combing \eqref{eq-It} and \eqref{eq-inter}, we get 
	\const{\CIt}
	\begin{align*}
		I(t)  \leq & C_{\CIt}\delta N^q I(t)+C\|f\|_{ \mL^{p}_{q}(T)}^{q}+C_\delta N^q \int_0^t I(\tau) \d \tau.
	\end{align*}
	Letting $\delta=\delta(N)$ be small enough so that $C_{\CIt} \delta N^q\leq 1/2$, we obtain that for all $t\in [0,T]$, 
	$$
	I(t)\leq C\|f\|_{ \mL^{p}_{q}(T)}^{q}+C \int_0^t I(\tau) \d \tau. 
	$$
	Grönwall's inequality yields 
	\begin{equation}\label{eq-W2-u}
	\|\p_t u\|_{\mL^{p}_q(T)} + \|u\|_{\mH^{2,p}_q(T)}\leq C I^{1/q}(T) \leq C\|f\|_{\mL^p_q(T)}. 
	\end{equation}
	
	(2). Let $v$ be the solution to 
	$$
	\p_t v= \tfrac{1}{2}\Delta v +f, \quad v(0)=0. 
	$$
	Again by \eqref{Eq-Heat-est}, one sees that 
	\begin{equation}\label{eq-W2-v}
	\|\p_t v\|_{\mH^{-1,p}_q(T)} + \|v\|_{\mH^{1,p}_q(T)} \leq C \|f\|_{\mH^{-1,p}_q(T)}. 
	\end{equation}
	Define $w:= u-v$. Then  
	$$
	\p_t w=\tfrac{1}{2} \Delta w+ b\cdot \nabla w+ b\cdot \nabla v, \quad w(0)=0. 
	$$ 
	Recalling that $p'=\frac{p_1p}{p_1+p}\in (1,p_1)$ and $q'=\frac{q_1q}{q_1+q}\in (1,q_1)$, by H\"older's inequality and \eqref{eq-W2-v}, we have 
	$$
	\|b\cdot\nabla v\|_{\mL^{p'}_{q'}(T)}\leq \|b\|_{\mL^{p_1}_{q_1}(T)} \|\nabla v\|_{\mL^p_q(T)}\leq C \|f\|_{\mH^{-1,p}_q(T)}. 
	$$
	This together with \eqref{eq-W2-u} implies  
	\begin{equation}\label{eq-W2-w}
	\|\p_t w\|_{\mL^{p'}_{q'}(T)} + \|w\|_{\mH^{2,p'}_{q'}(T)} \leq C\|f\|_{\mH^{-1,p}_q(T)}. 
	\end{equation}
	Using \eqref{Eq-PSobolev2} and noting that $\frac{d}{p'}+\frac{2}{q'}=\frac{d}{p}+\frac{d}{p_1}+\frac{2}{q}+\frac{2}{q_1}=1+\frac{d}{p}+\frac{2}{q}$,  one sees that 
	$$
	\|w\|_{\mH^{1,p}_q(T)} \overset{\eqref{Eq-PSobolev2}}{\leq}C \l( \|\p_t w\|_{\mL^{p'}_{q'}(T)}+ \|w\|_{\mH^{2,p'}_{q'}(T)}\r)\overset{\eqref{eq-W2-w}}{\leq }C \|f\|_{\mH^{-1,p}_{q}(T)}. 
	$$
	Combining \eqref{eq-W2-v} and the above estimate, we get $\|u\|_{\mH^{1,p}_{q}(T)}\leq \|v\|_{\mH^{1,p}_{q}(T)}+\|w\|_{\mH^{1,p}_{q}(T)}\leq C \|f\|_{\mH^{-1,p}_{q}(T)}$. So, we complete our proof. 
	\end{proof}
	
	\begin{remark}\label{Rek-local-W2}
	Let the assumptions in Theorem \ref{Th-Key1} or Theorem \ref{Th-Key2} hold. Suppose that $f\in \tH^{\a, p}_q(T)$ with $\a\in \{0,-1\}$. Then all the conclusions therein still hold if $\mH^{\cdots}_{\cdots}$ and $\mL^{\cdots}_{\cdots}$ are replaced by $\tH^{\cdots}_{\cdots}$ and $\tL^{\cdots}_{\cdots}$, respectively (cf. \cite{xia2020lqlp} or \cite{rockner2022weak}). 
	\end{remark}

	\section{Compactness criterion for \texorpdfstring{$L^2$}{} random fields}\label{Sec-Compact}
	In this section, we give a relative compactness criterion for the random fields on the Wiener-Sobolev space, which is essentially a consequence of \cite[Theorem 1 ]{bally2004relative}.  
	
	Let $(\Om,\sF, \bP)$ be a probability space. Assume $\{W_t\}_{t\in [0,T]}$ is a $d$-dimensional Brownian motion on $(\Om,\sF, \bP)$ and $\sF=\sigma\{W_t: t\in [0,T]\}$. $\bT=[0,T]\times \{1,2,\cdots,d\}$, $\mu$ is the product of the Lebesgue measure on $[0,T]$ times the uniform measure on $\{1,2,\cdots, d\}$.  $H:=L^2(\bT; \mu)$ and the scalar product is 
	$$
	\<f, g\>: = \sum_{i=1}^d \int_{0}^T f((t, i)) \, g((t,i)) \d t. 
	$$
	Let $I_m$ denote the multiple stochastic integral 
	\begin{equation}
		I_{m}\left(f_{m}\right)=m! \sum_{k_1,\cdots, k_m=1}^d \int\!\!\!\cdots\!\!\!\int_{0<t_1<\cdots<t_{m}<T} f_{m}\left((t_1, k_1), \cdots,(t_m, k_m)\right) \d W_{t_{1}}^{k_{1}} \ldots \d W_{t_{m}}^{k_{m}}
	\end{equation}
	of $L^2_s(\bT^m)$ (the collection of all symmetric elements in $L^2(\bT^m)$). Let $\sH_n$ denote the closed linear subspace of $L^2(\Om, \sF, \bP)$ generated by the random variables $\{H_n(I_1(h)): h\in H=L^2(\bT, \mu)\}$,  where $H_n$ is the $n$-th Hermite polynomial. The multiple integral $I_m$ is a map from $L^2_s(\bT^m)$ onto the Wiener chaos $\sH_m$ and any  $F\in L^2(\Om,\sF,\bP)$ can be expanded into a series of multiple stochastic integrals: $F= \sum_{m=0}^\infty I_m(f_m)$, where $I_0(F):= \bE F$. Let $\cS_p$ denote the class of smooth random variables $F=f(I_1(h_1),\cdots, I_1(h_m))$ and $f\in C_p^\infty(\R^d)$. The Malliavin derivative of a smooth random variable $F$ is the stochastic process $t\mapsto D_t F$ defined by
	$$
	D_t F:= \sum_{i=1}^m \p_{i} f(I_1(h_1),\cdots, I_1(h_m))h_i(t)
	$$
	Let $\mD^{1,2}$ be the closure of $\cS_p$ with respect to the norm
	$$
	\|F\|_{\mD^{1,2}}^2:= \bE F^2 + \bE \int_0^T |D_t F|^2 \d t. 
	$$
	Assume now $\cO$ is a bounded domain in $\R^d$ with smooth boundary. The Sobolev space $H^1_0(\mathcal{O})$ is defined as the closure of $C_c^\infty(\mathcal{O})$ with respect to the $H^1(\cO)$ norm given by
    $$
        \|f\|_{H^1(\mathcal{O})} = \l( \int_{\mathcal{O}} \left( |f(x)|^2 + |\nabla f(x)|^2 \right) \d x \r)^{\frac{1}{2}}.
    $$
    Let $F_n$ be a sequence of random fields in $L^2(\cO\times \Om)$.  
    The following result is a variant of a compactness criteria for sequences in $L^2(\cO\times \Om)$ due to Bally and Saussereau \cite{bally2004relative}.
	\begin{lemma}\label{Le-Comp}
	Assume $K>0$ and that the sequence $\{F_n\}_{n\in\mN}\subseteq  L^2 (\cO\times \Om)$ satisfies the following three conditions, for all $n\in \mN$: 
	\begin{align}\label{AS-1}
		\bE \|F_n\|_{H^1_x(\cO)}^2\leq K,  \tag{A$_1$}
	\end{align}
	\begin{align}\label{AS-2}
		\bE \int_{\cO}  \int_0^T   | D_sF_n(x) |^2 \d s\, \d x \leq K, \tag{A$_2$}
	\end{align}
	\begin{align}\label{AS-3}
		\bE \int_{\cO}\int_0^T\!\!\int_0^T \frac{| D_{s}  F_n(x)-D_{s'} F_n(x) |^2}{|s-s'|^{1+2\beta}} \d s\d s' \, \d x \leq K, \ \mbox{ for some } \beta>0, \tag{A$_3$}
	\end{align}
	then $\{F_n\}_{n\in \mN}$ is relatively compact in $L^2(\cO\times\Om)$. 
	\end{lemma}
\begin{proof}
    Since $\cO$ is a bounded smooth domain, there exists $\{e_k\}_{k\in \mN_+}\subseteq H_0^1(\cO)\cap C^\infty(\cO)$ and a sequence $\{\lambda_k\}_{k\in \mN_+}$ of positive real numbers with $\lambda_k\uparrow \infty \, (k\uparrow \infty)$, such that $\Delta e_k=-\lambda_k e_k$ and $\{e_k\}_{k\in \mN_+}$ forms a orthonormal base of $L^2(\cO)$ (cf. \cite{evans2010partial}). Moreover, $e_k/\sqrt{\lambda_k}$ forms a basis of $H^1_0(\cO)$ with norm $\|f\|_{H^1_0(\cO)}:= \l(\int_{\cO} |\nabla f|^2\r)^{1/2}$. Set $\<f, g\>:= \int_{\cO} fg$, then $F_n= \sum_{k=1}^\infty \<F_n, e_k\> e_k$. Integration by parts and \eqref{AS-1} yield, 
	\begin{equation*}
		\begin{aligned}
			&\l\| \sum_{{k=M}}^\infty \<F_n, e_k\> e_k \r\|_{L^2(\cO\times \Om)}^2=\bE \sum_{{k=M}}^\infty \<F_n, e_k\>^2 =\bE  \sum_{{k=M}}^\infty \lambda_k^{-2} \<F_n, \Delta e_k \>^2 \\
			=& \bE  \sum_{k= M}^\infty \lambda_k^{-1} \<\nabla F_{n}, \nabla e_k/\sqrt{\lambda_k}\>^{2}  \leq\lambda_M^{-1} \bE \sum_{k = M}^\infty   \< \nabla F_{n}, \nabla e_{k}/ \sqrt{\lambda_{k}}\>^2\\
			\leq & \lambda_M^{-1} \|\nabla F_n \|_{L^2(\cO \times \Om)}^2\leq C \lambda_M^{-1}\downarrow 0 \ \  (M\uparrow \infty). 
		\end{aligned}
	\end{equation*}
	Therefore, the relative compactness of the sequence $\{F_n\}_{n\in \mN}$ in $L^2(\cO\times \Om)$ reduces to the relative
	compactness of the sequence $\{\<F_n, e_k\>\}_{n\in \mN}$ in $L^2(\Om)$  for each $k\in \mN_+$. By \eqref{AS-1}, we have 
	\begin{equation}\label{eq-comp1}
	\bE \<F_n, e_k\>^2 \leq \bE \|F_n\|_{L^2(\cO)}^2 \leq K. 
	\end{equation}
	\eqref{AS-2} and \eqref{AS-3} yield for all $n\in\mN$ 
	\begin{equation}\label{eq-comp2}
	\begin{aligned}
		\bE \int_0^T | D_s \<F_n, e_k\> |^2  \d s=&\bE\int_0^T  \l| D_s \int_{\cO} F_n(x) e_k(x) \, \d x \r|^2  \d s\\
		\leq& \bE  \int_{\cO}\int_0^T |D_s F_n(x)|^2\, \d s\, \d x\leq K 
	\end{aligned}
	\end{equation}
	and 
	\begin{equation}\label{eq-comp3}
	\begin{aligned}
		&\bE\int_0^T\!\!\int_0^T \frac{|D_s\<F_n, e_k\>-D_{s'}\<F_n, e_k\>|^2}{|s-s'|^{1+2\beta}} \d s\d s'\\
		\leq &\bE \int_{\cO} \int_0^T\!\!\int_0^T \frac{|D_sF_n(x)-D_{s'}F_n(x)|^2}{|s-s'|^{1+2\beta}} \d s\d s'\,  \d x \leq K. 
	\end{aligned}
	\end{equation}
	By \eqref{eq-comp1}-\eqref{eq-comp3}, Theorem 1 and Lemma 1 of \cite{da1992compact} (with $\a\in (0, \beta\wedge \frac12)$ and $C=A_\alpha^{-1}$ therein), one sees that $\{\<F_n, e_k\>\}_{n\in \mN}$ is compact in $L^2(\Om)$ for each $k\in \mN_+$. So, we complete our proof. 
	\end{proof} 
	
	\section{Estimates for the case of regular coefficients}\label{Sec-Key}
	Throughout this section, we assume $b\in L^\infty([0,T]; C_b^2)$. The unique strong solution to SDE \eqref{Eq-SDE} with $s=0$ is denoted by $X_t^x$. Recall that $K_f(m)$, $K'_f(m)$ and $\om_f(\delta)$ are defined in \eqref{Eq-Kf}, \eqref{Eq-K'f} and \eqref{eq-om-f}, respectively. The main purpose of this section is to prove \begin{proposition}\label{Prop-Key}
	Let $d\geq 3$, $\{a(m)\}_{m\in \mN_+}$ be a sequence converging to zero and $\ell(\delta)$ be a monotonically increasing function on $(0,T)$ with $\lim_{\delta\downarrow0}\ell(\delta)= 0$. 
	\begin{enumerate}[(a)]
		\item Assume that $K_b(m)\leq a(m)$. Then for any $r\geq 2$, $p\in (\frac{d}{d-1}, d)$ and $\gamma\in (0,1/2)$, 
		\begin{equation}\label{Eq-nablaX}
		\l\| \nabla X_t^x-\mathrm{I} \r\|_{\widetilde L^{pr}_xL^r_\om} \leq C t^{\gamma/2r}, \quad \mbox{ for all } 0\leq t\leq T, 
		\end{equation}
		\begin{equation}\label{Eq-DX}
		\| D_{s} X_t^x-\mathrm{I} \|_{\widetilde L^{pr}_xL^r_\om} \leq C (t-s)^{\gamma/2r}, \quad \mbox{ for a.e. } s\in [0,T] \mbox{ with } 0\leq s\leq t\leq T  
		\end{equation}
		and 
		\begin{equation}\label{Eq-DX-Holder}
		\| D_{s} X_t^x -D_{s^{\prime}} X_t^x \|_{\widetilde L^{pr}_xL^r_\om} \leq C |s-s'|^{\gamma/4r}, \quad \mbox{ for a.e. } s, s' \in [0,T] \mbox{ with } 0\leq s, s'\leq t\leq T,  
		\end{equation}
		where $C$ only depends on $d, T, r,  p, \gamma, \|b\|_{\mL^d_\infty(T)}$ and $\{a(m)\}$.
		\item Assume that $p_1,q_1\in (2,\infty)$ and $d/p_1+2/q_1=1$, $K'_b(m)\leq a(m)$ and $\om_b(\delta)\leq \ell(\delta)$. Then for any $r\geq 2$, a.e. $s,s'\in [0,T]$ with $0\leq s, s'\leq t\leq T$, $p\in (\frac{p_1}{p_1-1}, p_1)$ and $\gamma\in (0,\frac{1}{2}-\frac{1}{q_1})$, the estimates \eqref{Eq-nablaX}-\eqref{Eq-DX-Holder} still hold, and the constant $C$ only depends on $d, p_1, q_1, T, r, p$, $\gamma, \{a(m)\}$ and $\ell(\delta)$. 
	\end{enumerate} 
	\end{proposition}
	\begin{remark}\label{Rek-Global}
	Using \eqref{Eq-Keyest1} and \eqref{Eq-Keyest2} below,  and following the proof for Proposition \ref{Prop-Key}, indeed one can verify that the $\widetilde L^{pr}_x L^r_\omega$-norm on the left-hand side of \eqref{Eq-nablaX}-\eqref{Eq-DX-Holder} can be replaced by the usual $ L^{pr}_x L^r_\omega$-norm if one further assume that $b\in \mL^p_{q}(T)$ with $q=\gamma^{-1}$ in case (a) and $q=(1/q_1+\gamma)^{-1}$ in case (b), and the constants on the right-hand side of \eqref{Eq-nablaX}-\eqref{Eq-DX-Holder} depend only $\|b\|_{\mL^p_q(T)}$ and parameters mentioned above. 
	\end{remark}
	The proof of Proposition \ref{Prop-Key} relies on the following lemma, which contains the key estimates of this paper. 
	\begin{lemma}\label{Le-Iteration}
	Let $d\geq 3$, $0\leq S_0\leq S_1\leq  T$ and $\{a(m)\}_{m\in \mN_+}$ be a sequence converging to zero. 
	\begin{enumerate}[(a)]
		\item 
		Suppose $b\in \mL^d_\infty(T)$ and $K_b(m)\leq a(m)$. Assume that $f_i\in  L^\infty([0,T]; C_b^2)\, (i\in \mN_+)$ 
		$$
		\sup_{i\in \mN_+} \|f_i\|_{\mL^d_\infty(T)}\leq N \ \mbox{ and }\   \sup_{t\in [0,T]; i\in \mN_+} K_{f_i}(m)  \leq d_m. 
		$$ 
		Then for any $p\in (\frac{d}{d-1},d)$, 
		$q\in (2,\infty)$, $\a_i\in \{1,2,\cdots d\}\, (i=1,2,\cdots)$  and all $n \in \mN_+$, it holds that 
		\begin{equation}\label{Eq-Keyest1}
		\begin{aligned}
			&\l\| \bE \int\!\!\!\cdots\!\!\!\int_{\Delta_n(S_0,S_1)} \prod_{i=1}^n \p_{\a_i} f_i\left(t_{i}, X_{t_{i}}^x\right) \ \d t_{1} \cdots \d t_{n} \r\|_{L^p_x}\\
			\leq & C^{n} \l(m^2N \sqrt{S_1-S_0}+d_m\r)^{n-1}  \|f_{n}\|_{\mL^{p}_q(S_0,S_1)}
		\end{aligned}
		\end{equation}
		and 
		\begin{equation}\label{Eq-Keyest1-1}
		\begin{aligned}
			&\l\| \bE \int\!\!\!\cdots\!\!\!\int_{\Delta_n(S_0,S_1)} \prod_{i=1}^n \p_{\a_i} f_i\left(t_{i}, X_{t_{i}}^x\right) \ \d t_{1} \cdots \d t_{n} \r\|_{\widetilde{L}^p_x}\\
			\leq & C^{n} \l(m^2N \sqrt{S_1-S_0}+d_m\r)^{n-1}  \|f_{n}\|_{\tL^{p}_q(S_0,S_1)},
		\end{aligned}
		\end{equation}
		where $C$ only depends on $d, p, q, T, \|b\|_{\mL^d_\infty}(T)$ and  $\{a(m)\}$. 
		
		\item  Suppose $b\in \mL^{p_1}_{q_1}(T)$ with $p_1,q_1\in (2,\infty)$ and $d/p_1+2/q_1=1$, and $K'_b(m)\leq a(m)$. Assume that $f_i\in \mL^{p_1}_{q_1}(T)\cap L^\infty([0,T]; C_b^2)\,  (i\in \mN_+)$. Then for any $p\in (\frac{p_1}{p_1-1}, p_1)$, $q\in (2,q_1)$, $\{\a_i\}_{i=1}^\infty\subseteq\{1,2,\cdots,d\}$ and all $n\in \mN_+$ 
		\begin{equation}\label{Eq-Keyest2}
		\begin{aligned}
			&\l\|  \bE \int\!\!\!\cdots\!\!\!\int_{\Delta_n(S_0,S_1)} \prod_{i=1}^n \p_{\a_i} f_i\left(t_{i}, X_{t_{i}}^x\right) \ \d t_{1} \cdots \d t_{n} \r\|_{L^p_x}\\
			\leq& C^{n+1} \l(\prod_{i=1}^{n-1}\|f_{i}\|_{\mL^{p_1}_{q_1}(S_0, S_1)} \r) \|f_n\|_{\mL^p_q(S_0,S_1)}
		\end{aligned}
		\end{equation} 
		and 
		\begin{equation}\label{Eq-Keyest2-1}
		\begin{aligned}
			&\l\|  \bE \int\!\!\!\cdots\!\!\!\int_{\Delta_n(S_0,S_1)} \prod_{i=1}^n \p_{\a_i} f_i\left(t_{i}, X_{t_{i}}^x\right) \ \d t_{1} \cdots \d t_{n} \r\|_{\widetilde{L}^p_x}\\
			\leq& C^{n+1} \l(\prod_{i=1}^{n-1}\|f_{i}\|_{\tL^{p_1}_{q_1}(S_0, S_1)} \r) \|f_n\|_{\tL^p_q(S_0,S_1)}, 
		\end{aligned}
		\end{equation}
		where $C$ only depends on $d, p_1, q_1, p, q, T$ and  $\{a(m)\}$, and when $n=1$, the right-hand side of the inequalities \eqref{Eq-Keyest2} and   \eqref{Eq-Keyest2-1} should be understood as $C^2 \|f_1\|_{\mL^p_q(S_0,S_1)}$ and $C^2 \|f_1\|_{\tL^p_q(S_0,S_1)}$, respectively 
	\end{enumerate}
	\end{lemma}
	\begin{proof} 
	We first reduce the left-hand side of \eqref{Eq-Keyest1} and \eqref{Eq-Keyest2} to some quantity associated with a family of parabolic equations. For fixed $n\in \mN_+$ and $\{\a_i\}_{i=1}^\infty \subseteq \{1,2,\cdots,d\}$, we set $u_{n+1}= 1$ and    for any $k\in \{1,2,\cdots, n\}$, let $g_{k}:= (\p_{\a_k} f_{k}) u_{k+1}$ and $u_k\in \cap_{p,q\in (1,\infty)}\tH^{2,p}_q(S_1)$ be the unique function solving equation 
	\begin{equation}\label{eq-uk}
	\p_t u_k + \tfrac{1}{2}\Delta u_k+ b\cdot \nabla u_k+ g_k=0 \ \mbox{in} \ (S_0,S_1)\times \R^d, \quad  u_k(S_1)=0
	\end{equation}
	(cf. \cite{xia2020lqlp}). Then the generalized It\^o formula yields 
	$$
	-u_k(t, X_t^x)= -\int_t^{S_1} g_k(s, X_s^x) \d s + \int_t^{S_1}\nabla u_k(s, X_s^x) \, \d W_s, \quad \forall t\in [0,S_1]
	$$
	which implies 
	\begin{equation}\label{eq-EFt-g}
	\bE^{\sF_t} \int_t^{S_1} g_k(s, X_s^x) \d s = u_k(t, X_t^x). 
	\end{equation}
	Here the conditional expectation $\bE(F|\sG)$ is denoted by $\bE^\sG F$. 
	By the Markov property and \eqref{eq-EFt-g}, 
	\begin{equation}\label{eq-Exp-iteration}
		\begin{aligned}
			& \bE^{\sF_{S_0}} \int\!\!\!\cdots\!\!\!\int_{\Delta_n(S_0,S_1)} \prod_{i=1}^n \p_{\a_i} f_i(t_i, X_{t_i}^x)  \ \dif t_1\cdots\dif t_{n}\\
			=&\bE^{\sF_{S_0}} \int\!\!\!\cdots\!\!\!\int_{\Delta_{n-1}(S_0,S_1)} \prod_{i=1}^{n-1} \p_{\a_i} f_i(t_i, X_{t_i}^x)\ \bE^{\sF_{t_{n-1}}}\left( \int^{S_1}_{t_{n-1}}\p_{\a_n}f_{n}(t_n, X_{t_{n}}^x)\dif t_n\right)\dif t_1\cdots\dif t_{n-1}\\
			\overset{\eqref{eq-EFt-g}}{=} &  \bE^{\sF_{S_0}} \int\!\!\!\cdots\!\!\!\int_{\Delta_{n-2}(S_0,S_1)} \prod_{i=1}^{n-2} \p_{\a_i} f_i(t_i, X_{t_i}^x) \\
            &\qquad \qquad \qquad \qquad \qquad\quad  \l[\int_{t_{n-2}}^{S_1}  (\p_{\a_{n-1}}f_{n-1}u_n)(t_{n-1}, X_{t_{n-1}}^x) \ \d t_{n-1} \r]  \d t_1\cdots\dif t_{n-2}\\
			=& \bE^{\sF_{S_0}} \int\!\!\!\cdots\!\!\!\int_{\Delta_{n-1}(S_0,S_1)} \prod_{i=1}^{n-2} \p_{\a_i} f_i(t_i, X_{t_i}^x) g_{n-1}(t_{n-1}, X_{t_{n-1}}^x) \ \d t_1\cdots \d t_{n-1}\\
			=&\cdots= u_1(S_0, X_{S_0}^x). 
		\end{aligned}
	\end{equation}
	We need to further simplify the expectation of the term on the right-hand side of the above equation.  Let $U$ be the solution to the following PDE: 
	\begin{equation}\label{eq-U}
	\p_t U=\tfrac{1}{2} \Delta U+B\cdot \nabla U+G \ \mbox{in} \ (0,S_1)\times \R^d, \quad U(0)=0, 
	\end{equation}
	where   
	$$
	B(t,x)=b(S_1-t,x)\1_{[0, S_1-S_0]}(t)+b(t+S_0-S_1,x)\1_{(S_1-S_0, S_1]}(t)
	$$
	and 
	\begin{equation}\label{eq-G}
	G(t,x)=g_1(S_1-t,x)\1_{[0, S_1-S_0]}(t).
	\end{equation}
	We note that $u_1(S_1-t)=U(t)$ for all $t\in [0, S_1-S_0]$ and that $V(t):=U(t+(S_1-S_0))$ satisfies 
	$$
	\p_t V=\tfrac{1}{2} \Delta V+b\cdot \nabla V \ \mbox{in} \ (0,S_0)\times \R^d, \quad V(0,x)=U(S_1-S_0,x)=u_1(S_0,x). 
	$$
	Therefore, for any $p\in [1,\infty)$, 
	\begin{equation}\label{eq-f-ite1}
	\begin{aligned}
		&\int_{\mR^d} \ \l| \bE \int\!\!\!\cdots\!\!\!\int_{\Delta_n(S_0,S_1)}\prod_{i=1}^n \p_{\a_i} f_i(t_i, X_{t_i}^x) \ \dif t_1\cdots\dif t_{n}  \r|^p \d x \\
		\overset{\eqref{eq-Exp-iteration}}{=}& \int_{A} |\bE u_1(S_0, X_{S_0}^x)|^p \d x =\int_{A} |V(S_0,x)|^p \d x=\|U(S_1)\|_{L^p_x}^p. 
	\end{aligned}
	\end{equation}
	Thanks to \eqref{eq-f-ite1}, now we can use the analytical results presented in Section 2 to prove the desired results, and we only give the proofs for \eqref{Eq-Keyest1} and \eqref{Eq-Keyest2} as the proofs for \eqref{Eq-Keyest1-1} and \eqref{Eq-Keyest2-1} are almost the same. 
	
	\noindent {\bf Case (a): $b\in C([0,T];L^d)$.} 
	Set 
	$$
	f_{k, m}(t):= f_k(t)*\rho_m, \quad \bar f_{k,m}:= f_k-f_{k,m}. 
	$$
	Let 
	$$p\in (d/(d-1),d)\ \mbox{ and } \ q \in (2,\infty).
	$$ 
	By the definitions of $g_k$ and $u_{k+1}$, 
	\begin{equation}\label{eq-gk-gk+1}
	\begin{aligned}
		&\|g_{k}\|_{\mH^{-1,p}_q(S_0,S_1)}=\|(\p_{\a_k}f_k) u_{k+1}\|_{\mH^{-1,p}_q(S_0,S_1)} \\
		\leq &
		\| (\p_{\a_k}f_{k, m})\, u_{k+1}\|_{\mH^{-1,p}_q(S_0,S_1)} + \| (\p_{\a_k} \bar f_{k, m})\, u_{k+1}\|_{\mH^{-1,p}_q(S_0,S_1)}=:I_1+I_2. 
	\end{aligned}
	\end{equation}
	Recalling that $u_{k+1}$ solves \eqref{eq-uk} with $k$ replaced by $k+1$, using \eqref{Eq-W2}, we get 
	$$
	(S_1-S_0)^{-1}\|u_{k+1}\|_{\mH^{-1, p}_q(S_0, S_1)}+\|u_{k+1}\|_{\mH^{1, p}_q(S_0, S_1)} \leq C \|g_{k+1}\|_{\mH^{-1, p}_q(S_0, S_1)}. 
	$$
	An interpolation inequality yields 
	\begin{align*}
		\|u_{k+1}\|_{\mL^{p}_q(S_0, S_1)}\leq&  C \|u_{k+1}\|_{\mH^{-1, p}_q(S_0, S_1)}^{1/2} \|u_{k+1}\|_{\mH^{1,p}_q(S_0, S_1)}^{1/2}\\
		\leq& C \sqrt{S_1-S_0} \|g_{k+1}\|_{\mH^{-1, p}_q(S_0, S_1)}. 
	\end{align*}
	Thus, like in the proof for \eqref{eq-bmu'}, we have 
	\begin{equation}\label{eq-fmu-I1}
	\begin{aligned}
		I_1\leq& C\| \p_{\a_k}f_{k, m}\, u_{k+1}\|_{\mL^{p}_q(S_0,S_1)} \\
		\leq & C \|f_{k, m}\|_{L^\infty([0,T]; C_b^1)}\|u_{k+1} \|_{\mL^{p}_q(S_0,S_1)} \\\leq& C m^2 N  \sqrt{S_1-S_0}\|g_{k+1}\|_{\mH^{-1,p}_q(S_0,S_1)}. 
	\end{aligned} 
	\end{equation}
	On the other hand, recalling that $p\in (d/(d-1), d)$ (this implies $pd/(d+p)>1$), by the Sobolev embedding, 
	\begin{equation}\label{eq-fmu-I2}
	\begin{aligned}
		I_2\leq& \|\p_{\a_k} (\bar f_{k,m} u_{k+1})\|_{\mH^{-1,p}_q(S_0,S_1)}+ \|\bar f_{k,m}\p_{\a_k} u_{k+1}\|_{\mL^{\frac{dp}{d+p}}_q(S_0,S_1)}\\
		\leq &C \|\bar f_{k, m}\|_{\mL^{d}_\infty(S_0,S_1)} \|u_{k+1}\|_{\mL^{\frac{pd}{d-p}}_q(S_0,S_1)}+C\|\bar f_{k,m}\|_{\mL^d_\infty(T)} \|u_{k+1}\|_{\mH^{1,p}_q(S_0,S_1)} \\
		\leq& C d_m \|g_{k+1}\|_{\mH^{-1, p}_q(S_0,S_1)}. 
	\end{aligned}
	\end{equation}
	Combing \eqref{eq-gk-gk+1}-\eqref{eq-fmu-I2}, we get 
	\begin{align*}
		\|g_{k}\|_{\mH^{-1,p}_q(S_0,S_1)}\leq C \l(m^2N \sqrt{S_1-S_0}+d_m\r)\|g_{k+1}\|_{\mH^{-1,p}_q(S_0,S_1)},
	\end{align*}
	where $C$ only depends on $d, p, \gamma, T, \|b\|_{\mL^d_\infty(T)}$  and $\{a(m)\}$.  Recalling that $G$ is defined in \eqref{eq-G}, by the above estimate we obtain 
	\begin{equation}\label{eq-G-est1}
	\begin{aligned}
		&\|G\|_{\mH^{-1,p}_{q}(S_1)} \leq \|g_1\|_{\mH^{-1,p}_{q}(S_0, S_1)} \\
		\leq& C^{n}\l(m^2N \sqrt{S_1-S_0}+d_m\r)^{n-1} \|g_{n}\|_{\mH^{-1,p}_q(S_0,S_1)}\\
		=& C^{n}\l(m^2N \sqrt{S_1-S_0}+d_m\r)^{n-1} \|f_{n}\|_{\mL^{p}_q(S_0,S_1)}. 
	\end{aligned}
	\end{equation}
	Thus, 
	\begin{equation*}
		\begin{aligned}
			&\l\| \bE \int\!\!\!\cdots\!\!\!\int_{\Delta_n(S_0,S_1)} \prod_{i=1}^n \p_{\a_i} f_i(t_i, X_{t_i}^x)  \ \dif t_1\cdots\dif t_{n} \r\|_{{L}^p_x}\\
			\overset{\eqref{eq-f-ite1}}{\leq}& C \| U (S_1)\|_{{L}^p}
			\overset{\eqref{Eq-PMorrey}}{\leq}  C \l(\|\p_t U\|_{\mH^{-1,p}_q(S_1)}+\|U\|_{\mH^{1, p}_q(S_1)}\r)\\
			&(\mbox{ taking } \a=-1, \theta=1/2, \mbox{ and noticing } q>2)\\
			\overset{\eqref{Eq-W2}}{\leq} &C\|G\|_{\mH^{-1,p}_q(S_1)} \overset{\eqref{eq-G-est1}}{\leq}  C^{n} \l(m^2N \sqrt{S_1-S_0}+d_m\r)^{n-1}  \|f_{n}\|_{\mL^{p}_q(S_0,S_1)}, 
		\end{aligned}
	\end{equation*}
	where $C$ only depends on $d, p, \gamma,  \|b\|_{\mL^d_\infty}(T)$ and $\{a(m)\}$. So,  we complete the proof for \eqref{Eq-Keyest1}. 
	\
	\\
	{\bf Case (b): $b\in \mL^{p_1}_{q_1}(T)$.} Set  
	\begin{equation}\label{eq:pq}
	    p\in \l(\tfrac{p_1}{p_1-1}, p_1\r), \quad 
	q \in (2, q_1), \quad p'=\frac{p_1p}{p_1+p}, \quad q'=\frac{q_1q}{q_1+q}.
	\end{equation}
	Noting that $p_1\in (d,\infty)$ and $q_1\in (2,\infty)$, one sees that $p\in (1, p_1), q\in (2,q_1)$ and $p'\in (1,p), q'\in (1,q)$.  
	We claim that for each $k\in \{n, \cdots, 1\}$, $g_k$ can be written as the sum of two functions $g'_k$ and $g''_k$, and 
	\begin{equation}\label{eq-gk}
	\|g_k'\|_{\mH^{-1,p}_q(S_0, S_1)}+\|g_k''\|_{\mL^{p'}_{q'}(S_0, S_1)}  \leq C^{n-k+1} \l( \prod_{i=k}^{n-1}\|f_{i}\|_{\mL^{p_1}_{q_1}(S_0, S_1)}\r) \|f_n\|_{\mL^{p}_{q}(S_0, S_1)}, 
	\end{equation}
	where $C$ does not depend on $n$, and the right-hand side of the above inequality should be understood as $C\|f_n\|_{\mL^p_q(S_0,S_1)}$ when $k=n$. Recalling that $g_n=\p_{\a_n} f_n$, by letting $g_n'=g_n=\p_{\alpha_n} f_n$, one sees that 
	$$
	\|g'_n\|_{\mH^{-1,p}_q(S_0, S_1)}\leq C \|f_n\|_{\mL^{p}_{q}(S_0, S_1)}, \quad g''_n=0. 
	$$ 
	Thus, \eqref{eq-gk} holds for the case that $k=n$.  
	Suppose $n\geq 2$. Assume $g_k=g_k'+g_k''$ and that \eqref{eq-gk} holds for some $k\in \{n,\cdots, 2\}$. Then  $u_k$ can be decomposed as $u_k=u_k'+u_k''$, where $u_k'$ and $u_k''$ solve \eqref{eq-uk} with $g_k$ replaced by $g_k'$ and $g_k''$, respectively. By Theorem \ref{Th-Key2}, one sees that $u_k$ can be further decomposed as
	$$
	u_k=u_k'+u_k''=v_k'+w_k'+u_k'', 
	$$ 
	where $v_k', w_k'$ and $u_k''$ satisfy 
	\begin{equation}\label{eq-vk'}
	\|\p_t v_k'\|_{\mH^{-1,p}_q(S_0, S_1)} + \|v_k'\|_{\mH^{1,p}_q(S_0, S_1)} \leq C \|g_k'\|_{\mH^{-1, p}_q(S_0, S_1)}, 
	\end{equation}
	\begin{equation}\label{eq-wk'}
	\|\p_t w_k'\|_{\mL^{p'}_{q'}(S_0, S_1)} + \|w_k'\|_{\mH^{2,p'}_{q'}(S_0, S_1)} \leq C\|g_k'\|_{\mH^{-1, p}_q(S_0, S_1)}   
	\end{equation}
	and 
	\begin{equation}\label{eq-uk''}
	\|\p_t u_k''\|_{\mL^{p'}_{q'}(S_0, S_1)} + \|u_k''\|_{\mH^{2,p'}_{q'}(S_0, S_1)} \leq C\|g_k''\|_{\mL^{p'}_{q'}(S_0, S_1)}. 
	\end{equation}
	Let $r=\frac{p_1p}{p_1-p}$ and $s=\frac{1}{\gamma}=\frac{q_1q}{q-q_1}$. Recalling that $2<\frac{d}{p'}+\frac{2}{q'}=1+\frac{d}{p}+\frac{2}{q}=2+\frac{d}{r}+\frac{2}{s}$, due to \eqref{Eq-PSobolev1}, \eqref{Eq-PSobolev2}, we have 
	\begin{equation}\label{eq-uk-est1}
	\begin{aligned}
		&\|u_k\|_{\mL^r_s(S_0, S_1)}\leq \|v_k'\|_{\mL^r_s(S_0, S_1)}+\|w_k'\|_{\mL^r_s(S_0, S_1)}+\|u_k''\|_{\mL^r_s(S_0, S_1)}\\
		\overset{\eqref{Eq-PSobolev1}, \eqref{Eq-PSobolev2}}{\leq}& C \l(\|\p_t v_k'\|_{\mH^{-1,p}_q(S_0, S_1)} + \|v_k'\|_{\mH^{1,p}_q(S_0, S_1)} \r)+C \l(\|\p_t w_k'\|_{\mL^{p'}_{q'}(S_0, S_1)} + \|w_k'\|_{\mH^{2,p'}_{q'}(S_0, S_1)} \r)\\
		&+C \l(\|\p_t u_k''\|_{\mL^{p'}_{q'}(S_0, S_1)} + \|u_k''\|_{\mH^{2,p'}_{q'}(S_0, S_1)} \r) \\
		\overset{\eqref{eq-vk'}-\eqref{eq-uk''}}{\leq}& C \l( \|g_k'\|_{\mH^{-1, p_1}_{q_1}(S_0, S_1)}+ \|g_k''\|_{\mL^{p_1}_{q_1}(S_0, S_1)}\r)
	\end{aligned}
	\end{equation}
	and 
	\begin{equation}\label{eq-uk-est2}
	\|u_k\|_{\mH^{1, p}_q(S_0, S_1)} \leq C \l( \|g_k'\|_{\mH^{-1, p_1}_{q_1}(S_0, S_1)}+ \|g_k''\|_{\mL^{p_1}_{q_1}(S_0, S_1)}\r). 
	\end{equation}
	Set $g_{k-1}'=\p_{\a_{k-1}}(f_{k-1} u_k)$ and 
	$g_{k-1}''=-f_{k-1}(\p_{\a_{k-1}} u_k)$. By H\"older's inequality, \eqref{eq-uk-est1} and \eqref{eq-uk-est2}, we get 
	\begin{align*}
		&\|g_{k-1}'\|_{\mH^{-1,p}_q(S_0, S_1)} + \|g_{k-1}''\|_{\mL^{p'}_{q'}(S_0, S_1)}\\
		\leq& \|f_{k-1}\|_{\mL^{p_1}_{q_1}(S_0, S_1)} \l( \|u_k\|_{\mL^r_s(S_0, S_1)}+\|\nabla u_k\|_{\mL^p_q(S_0, S_1)} \r)\\
		\overset{\eqref{eq-uk-est1}, \eqref{eq-uk-est2}}{\leq}& C \|f_{k-1}\|_{\mL^{p_1}_{q_1}(S_0, S_1)} \l( \|g_k'\|_{\mH^{-1, p_1}_{q_1}(S_0, S_1)}+ \|g_k''\|_{\mL^{p_1}_{q_1}(S_0, S_1)}\r)\\
		\leq &  C^{n-k+2}\prod_{i=k-1}^{n-1}\|f_{i}\|_{\mL^{p_1}_{q_1}(S_0, S_1)}\|f_n\|_{\mL^q_q(S_0,S_1)}. 
	\end{align*}
	So, by induction \eqref{eq-gk} holds for all $k\in \{1, 2, \cdots, n\}$. In particular, 
	$$
	\|g_{1}'\|_{\mH^{-1,p}_q(S_0, S_1)} + \|g_{1}''\|_{\mL^{p'}_{q'}(S_0, S_1)} \leq C^{n+1}\l( \prod_{i=1}^{n-1}\|f_{i}\|_{\mL^{p_1}_{q_1}(S_0, S_1)}\r)\|f_n\|_{\mL^p_q(S_0,S_1)}. 
	$$
	Recalling that $G$ is defined in \eqref{eq-G}, this can be written as 
	\begin{align*}
		G(t,x)=&G'(t,x)+G''(t,x)\\
		:=& g_1'(S_1-t,x)\1_{[0, S_1-S_0]}(t)+g_1''(S_1-t,x)\1_{[0, S_1-S_0]}(t) 
	\end{align*}
	and 
	\begin{equation}\label{eq-EstG}
	\begin{aligned}
		\|G'\|_{\mH^{-1,p}_{q}(S_1)}+ \|G''\|_{\mL^{p'}_{q'}(S_1)}  \leq&  C\l( \|g_1'\|_{\mH^{-1,p}_{q}(S_0, S_1)}+ \|g_1''\|_{\mL^{p'}_{q'}(S_0, S_1)} \r) \\
		\leq& C^{n+1}\l(\prod_{i=1}^{n-1}\|f_{i}\|_{\mL^{p_1}_{q_1}(S_0, S_1)} \r)\|f_n\|_{\mL^p_q(S_0,S_1)}. 
	\end{aligned}
	\end{equation}
	Assume $U'$ and $U''$ solve \eqref{eq-U} with $G$ replaced by $G'$ and $G''$, respectively. As in the above argument, we see that $U=U'+U''=V'+W'+U''$ and that 
	\begin{equation}\label{eq-VWU}
	\begin{aligned}
		&\l(\|\p_t V'\|_{\mH^{-1,p}_q(S_1)} + \|V'\|_{\mH^{1,p}_q(S_1)}\r)+  \l(
		\|\p_t W'\|_{\mL^{p'}_{q'}(S_1)} + \|W'\|_{\mH^{2,p'}_{q'}(S_1)} \r)\\
		&\qquad \qquad \qquad \qquad \qquad \qquad  \quad+ \l(\|\p_t U''\|_{\mL^{p'}_{q'}(S_1)} + \|U''\|_{\mH^{2,p'}_{q'}(S_1)} \r)\\
		\leq& C\l( \|G'\|_{\mH^{-1,p}_{q}(S_1)}+ \|G''\|_{\mL^{p'}_{q'}(S_1)} \r) \\
		\overset{\eqref{eq-EstG}}{\leq}& C^{n+1} \l(\prod_{i=1}^{n-1}\|f_{i}\|_{\mL^{p_1}_{q_1}(S_0, S_1)} \r)\|f_n\|_{\mL^p_q(S_0,S_1)}, 
	\end{aligned}
	\end{equation}
	where the first inequality is due to Theorem \ref{Th-Key2}. 
	
	Recalling that $q>2$, by taking $\a=-1$ and $\theta=1/2$ in \eqref{Eq-PMorrey}, we get 
	\begin{equation}\label{eq-V'}
	\|V'\|_{\mL^p_\infty(S_1)} \leq   C \l(\|\p_t V'\|_{\mH^{-1,p}_q(S_1)}+ \| V'\|_{\mH^{1,p}_q(S_1)}\r). 
	\end{equation}
	Similarly, using inequality \eqref{Eq-PMorrey} with \(\alpha = 0\), \(\theta = \frac{d}{2p_1}\), and \(q\) therein taken as \(q'\) (given by \eqref{eq:pq}), and noting that $1-\frac{1}{q'}=1-\frac{1}{q_1}-\frac{1}{q}=\frac{d}{2p_1}+\frac{1}{2}-\frac{1}{q}>\frac{d}{2p_1}=\theta$ and $\frac{1}{p}=\frac{1}{p'}-\frac{2\theta}{d}$, we get 
	\begin{equation}\label{eq-WU}
	\begin{aligned}
		&\|W'\|_{\mL^{p}_\infty(S_1)}+\|U''\|_{\mL^{p}_\infty(S_1)} \leq C \l( \|W'\|_{\mH^{2\theta, p'}_\infty(S_1)}+\|U''\|_{\mH^{2\theta, p'}_\infty(S_1)} \r) \\
		\leq&  C \l(\|\p_t W'\|_{\mL^{p'}_{q'}(S_1)}+\| W'\|_{\mH^{2,p'}_{q'}(S_1)}\r)+  C \l(\|\p_t U''\|_{\mL^{p'}_{q'}(S_1)}+ \| U''\|_{\mH^{2,p'}_{q'}(S_1)}\r). 
	\end{aligned}
	\end{equation}
	Combining \eqref{eq-VWU}-\eqref{eq-WU} with \eqref{eq-f-ite1}, we obtain 
	\begin{align*}
		&\l\| \bE \int\!\!\!\cdots\!\!\!\int_{\Delta_n(S_0,S_1)}\prod_{i=1}^n \p_{\a_i} f_i(t_i, X_{t_i}^x) \ \dif t_1\cdots\dif t_{n}  \r\|_{L^p_x}\\
		\leq& C \|U\|_{\mL^{p}_\infty(S_1)} \leq C\l( \|V'\|_{\mL^{p}_\infty(S_1)}+  \|W'\|_{\mL^{p}_\infty(S_1)}+  \|U''\|_{\mL^{p}_\infty(S_1)}\r) \\
		\leq& C^{n+1}\l(\prod_{i=1}^{n-1}\|f_{i}\|_{\mL^{p_1}_{q_1}(S_0, S_1)} \r)\|f_n\|_{\mL^p_q(S_0,S_1)}. 
	\end{align*}
	So, we complete our proof. 
	\end{proof}
	
	\medskip
	
	Note that $b\in L^\infty([0,T]; C_b^2)$, the solution to the SDE \eqref{Eq-SDE} is differentiable with respect to $x$, and $\nabla X_t^x$ satisfies 
	$$
	\nabla X_t^x= {\rm I}+ \int_0^t \nabla b(s, X_s^x) \nabla X_s^x\  \d s. 
	$$
	Regarding the above equation as a linear random ODE for $\nabla X_t^x$, this 
	equation has a unique solution and it is given by 
	\begin{equation}\label{Eq-X'exp}
	\nabla X_t^x={\rm I}+\sum_{n=1}^{\infty} \int\!\!\!\cdots\!\!\!\int_{\Delta_n(t)} \nabla b\left(t_{n}, X_{t_{n}}^x\right)\cdots\nabla b\left(t_{1}, X_{t_{1}}^x\right)\,   \d t_1 \cdots \d t_n, 
	\end{equation}
	provided that this series is convergent (cf. \cite{mohammed2015sobolev}).  Moreover, for any $0\leq t_0\leq t\leq T$,  
	\begin{equation}\label{Eq-X'exp2}
		\nabla X_t^x=\nabla X_{t_0}^x+\sum_{n=1}^{\infty} \int\!\!\!\cdots\!\!\!\int_{\Delta_n(t_0, t)} \nabla b\left(t_{n}, X_{t_{n}}^x\right)\cdots\nabla b\left(t_{1}, X_{t_{1}}^x\right) \, \nabla X_{t_0}^x \ \d t_1 \cdots \d t_n. 
	\end{equation}
	On the other hand, the Malliavin derivative $D_s X_t^x$  is the solution of the linear stochastic equation
	$$
	D_{s} X_t^x={\rm I}+\int_{s}^{t} \nabla b \left(r, X_r^x\right) D_{s} X_r^x\ \d r, 
	$$
	for a.e. $s\in [0,T]$ with $s\leq t$, and $D_sX_t=0$ for a.e. $s\in [0,T]$ with $s>t$. Thus,  one sees that
	\begin{equation}\label{Eq-DXexp}
	D_s X_t^x= {\rm I}+ \sum_{n=1}^{\infty} \int\!\!\!\cdots\!\!\!\int_{\Delta_n(s, t)}   
	\nabla b\left(t_{n}, X_{t_{n}}^x\right)\cdots\nabla b\left(t_{1}, X_{t_{1}}^x\right)\ \d t_1 \cdots \d t_n,
	\end{equation}
	for a.e. $s\in [0,T]$ with $s\leq t$, and
	\begin{equation*}
		\begin{aligned}
			D_{s} X_t^x-D_{s^{\prime}} X_t^x =&\int_{s}^{t} \nabla b\left(r, X_r^x\right) D_{s} X_r^x \ \d r-\int_{s^{\prime}}^{t} \nabla b\left(r, X_r^x\right) D_{s^{\prime}} X_r^x \ \d r \\
			=&\int_{s}^{s^{\prime}} \nabla b\left(r, X_r^x\right) D_{s} X_r^x \ \d r+\int_{s^{\prime}}^{t} \nabla b\left(r, X_r^x\right)\left(D_{s} X_r^x-D_{s^{\prime}} X_r^x\right)\  \d r \\
			=&D_{s} X_{s^{\prime}}^x-{\rm I}+\int_{s^{\prime}}^{t} \nabla b\left(r, X_r^x\right)\left(D_{s} X_r^x-D_{s^{\prime}} X_r^x\right) \ \d r
		\end{aligned}
	\end{equation*}
	for a.e. $s,s'\in [0,T]$ with $s<s'\leq t$. Iterating, we get 
	\begin{equation}\label{Eq-DXexp2}
		\begin{aligned}
			D_{s} X_t^x -D_{s^{\prime}} X_t^x =&\left({\rm I}+\sum_{n=1}^{\infty} \int\!\!\!\cdots\!\!\!\int_{\Delta_n(s^{\prime}, t)} \nabla b\left(t_{n}, X_{t_{n}}^x\right)\cdots\nabla b\left(t_{1}, X_{t_{1}}^x\right)\ \d t_1 \cdots \d t_n\right) \cdot \left(D_{s} X_{s^{\prime}}^x-{\rm I}\right)\\
			\overset{\eqref{Eq-DXexp}}{=}& D_{s'}X_t^x\cdot \left(D_{s} X_{s^{\prime}}^x-{\rm I}\right), 
		\end{aligned}
	\end{equation}
	for a.e. $s,s'\in [0,T]$ with $s<s'\leq t$.  
	
	We are now in a position to prove our Proposition 
	\ref{Prop-Key}. 
	\begin{proof}[Proof of Proposition \ref{Prop-Key}]
	{\bf Case (a).}  We only need to prove the case where $r$ is a positive even integer. For any $n\in \mN_+$ and  $0\leq S_0\leq S_1 \leq  T$, it is not hard to see that  
	$$
	\l( \int\!\!\!\cdots\!\!\!\int_{\Delta_n(S_0, S_1)}  \p_{\a_{n-1}} b^i(t_n, X_{t_n}^x)\cdot \p_{\a_{n-2}} b^{\a_{n-1}}(t_{n-1}, X_{t_{n-1}}^x)\cdots \p_j b^{\a_{1}}(t_1, X_{t_{1}}^x)\  \d t_1\d t_2 \cdots \d t_{n} \r)^r
	$$
	can be written as a sum of at most ${r^{rn}}$ terms of the form 
    $$
        \int\!\!\!\cdots\!\!\!\int_{\Delta_{rn}(S_0, S_1)}  \p_{\beta'_1} b^{\beta_1}(t_1, X_{t_1}^x)\cdot \p_{\beta'_2} b^{\beta_2}(t_2, X_{t_2}^x)\cdots \p_{\beta'_{rn}} b^{\beta_{rn}}(t_{rn}, X_{t_{rn}}^x) \ \d t_1\d t_2\cdots \d t_{rn}.   
    $$
    
    This is because 
    there are \(\binom{rn}{n} \binom{(r-1)n}{n}\cdots \binom{n}{n}=\frac{(rn)!}{(n!)^r} (\leq r^{rn})\) ways to choose \(n\) elements at a time from \(r\times n\) elements, repeated \(r\) times.
    
    Fix $p\in (d/(d-1),d)$ and $q=\gamma^{-1}\in (2,\infty)$. By the above discussion and \eqref{Eq-Keyest1}, we have
	\const{\CXiteone}
	\begin{equation*}
		\begin{aligned}
			&\l\|\int\!\!\!\cdots\!\!\!\int_{\Delta_n(S_0, S_1)} \nabla b\left(t_{n}, X_{t_{n}}^x\right)\cdots\nabla b\left(t_{1}, X_{t_{1}}^x\right)\ \d t_1 \cdots \d t_n \r\|_{L^{pr}_x(B_1(z))L^r_\om}\\
			\leq &C \sum_{i, j=1}^{d} \sum_{\a_1, \cdots, \a_{n-1}=1}^{d} \Big\| \int\!\!\!\cdots\!\!\!\int_{\Delta_n(S_0, S_1)}  \p_{\a_{n-1}} b^i(t_n, X_{t_n}^x)\cdot \p_{\a_{n-2}} b^{\a_{n-1}}(t_{n-1}, X_{t_{n-1}}^x)\cdots   \\
			& \qquad\qquad\qquad\qquad \qquad\qquad\qquad \qquad \cdot \p_j b^{\a_{1}}(t_1, X_{t_{1}}^x)\  \d t_1\d t_2 \cdots \d t_{n}   \Big\|_{L^{pr}_x(B_1(z))L^r_\om}\\
			=& C \sum_{i, j=1}^{d} \sum_{\a_1, \cdots, \a_{n-1}=1}^{d} \Big[ \int_{B_1(z)} \Big( \sum_{\beta, \beta'}\bE \int\!\!\!\cdots\!\!\!\int_{\Delta_{rn}(S_0, S_1)}  \p_{\beta'_1} b^{\beta_1}(t_1, X_{t_1}^x)\cdot \p_{\beta'_2} b^{\beta_2}(t_2, X_{t_2}^x)\cdots \\
			&\qquad \qquad \qquad\qquad \qquad \qquad \qquad \qquad \cdot \p_{\beta'_{rn}} b^{\beta_{rn}}(t_{rn}, X_{t_{rn}}^x) \d t_1\d t_2\cdots \d t_{rn} \Big)^p \d x\Big]^{1/pr}\\
			\leq &  C \sum_{i, j=1}^{d} \sum_{\a_1, \cdots, \a_{n-1}=1}^{d}  \Big[\sum_{\beta, \beta'} \Big\| \bE \int\!\!\!\cdots\!\!\!\int_{\Delta_{rn}(S_0, S_1)}  \p_{\beta'_1} b^{\beta_1}(t_1, X_{t_1}^x)\cdot \p_{\beta'_2} b^{\beta_2}(t_2, X_{t_2}^x)\cdots \\
			& \qquad \qquad\qquad \qquad \qquad \qquad \qquad\qquad  \cdot \p_{\beta'_{rn}} b^{\beta_{rn}}(t_{rn}, X_{t_{rn}}^x) \d t_1\d t_2\cdots \d t_{rn} \Big\|_{L^p_x(B_1(z))}  \Big]^{1/r}\\
			\overset{\eqref{Eq-Keyest1-1}}{\leq}& (rC_{\CXiteone})^{n} \l(m^2 \|b\|_{\mL^d_\infty(T)}\sqrt{S_1-S_0}+a(m)\r)^{n-1/r}  \|b\|_{\mL^d_\infty(T)}^{1/r} (S_1-S_0)^{\gamma/r}.  
		\end{aligned}
	\end{equation*}
	Here we also used the fact that the sum $\sum_{\beta, \beta'}$ contains at most ${r^{rn}}$ terms. The constant $C_{\CXiteone}>1$ only depends on $d, T, \|b\|_{\mL^d_\infty(T)}, \{a(m)\}, r, p$ and $\gamma$. Letting $m$ be large enough such that $C_{\CXiteone}r a(m)\leq 1/4$ and then choosing $T_r=(4C_{\CXiteone} r m^2 \|b\|_{\mL^d_\infty(T)})^{-2}>0$ such that $C_{\CXiteone} r m^2 \|b\|_{\mL^d_\infty(T)}\sqrt{T_r}= 1/4$, we have for any  $0\leq S_1-S_0\leq T_r$, 
	\begin{equation}\label{eq-prod-b'} 
	\begin{aligned} 
		&\sum_{n=1}^{\infty} \l\|\int\!\!\!\cdots\!\!\!\int_{\Delta_n(S_0, S_1)} \nabla b\left(t_{n}, X_{t_{n}}^x\right)\cdots\nabla b\left(t_{1}, X_{t_{1}}^x\right)\ \d t_{1} \cdots \d t_{n} \r\|_{\widetilde{L}^{pr}_xL^r_\om}\\
		\leq &  \sum_{n=1}^{\infty} (C_{\CXiteone} r)^n \l(m^2 \|b\|_{\mL^d_\infty(T)} \sqrt{T_r}+a(m)\r)^{n-1/r}  \|b\|_{\mL^d_\infty(T)}^{1/r} (S_1-S_0)^{\gamma/r}\\
		\leq & 2 (rC_{\CXiteone})^{1/r} \|b\|_{\mL^d_\infty(T)}^{1/r} (S_1-S_0)^{\gamma/r} .
	\end{aligned}
	\end{equation}
	  Using \eqref{Eq-X'exp} and applying \eqref{eq-prod-b'} with \(S_0=0\) and \(S_1=t\), for each $t\in [0,T_r]$,
    \begin{equation}\label{Eq-X'-Lr}
	\begin{aligned}
        \l\| \nabla X_t^x-\mathrm{I} \r\|_{\widetilde{L}^{pr}_xL^r_\om} \overset{\eqref{Eq-X'exp}}{\leq}&  C \sum_{n=1}^{\infty}  \l\|\int\!\!\!\cdots\!\!\!\int_{\Delta_n(t)} \nabla b\left(t_{n}, X_{t_{n}}^x\right)\cdots\nabla b\left(t_{1}, X_{t_{1}}^x\right) \  \d t_{1} \cdots \d t_{n} \r\|_{\widetilde{L}^{pr}_xL^r_\om} \\
	\overset{\eqref{eq-prod-b'}}{\leq}& C t^{\gamma/r}. 
	\end{aligned}
    \end{equation}
    In particular, we have \(\l\| \nabla X_t^x-\mathrm{I} \r\|_{\widetilde{L}^{pr}_xL^r_\om} \leq C t^{\gamma/2r}\), \(t\in [0,T_{2r}]\). We then extend this result to the entire interval \([0,T]\), showing that
    \begin{equation}\label{eq-DX-Lr}
        \l\| \nabla X_t^x-\mathrm{I} \r\|_{\widetilde{L}^{pr}_xL^r_\om}\leq  C t^{\gamma/2r}, \quad t\in [0,T].  
    \end{equation}
    For any $t \in [T_{2r}, 2T_{2r}\wedge T]$, by \eqref{Eq-X'exp2}, H\"older's inequality, \eqref{eq-prod-b'} and \eqref{Eq-X'-Lr} we get 
	\begin{align*}
		&\l\| \nabla X_t^x -\mathrm{I}\r\|_{\widetilde{L}^{pr}_xL^r_\om} \\
		\leq&  \l\| \nabla X_{T_{2r}}^x-\mathrm{I} \r\|_{{\widetilde{L}^{pr}_xL^r_\om}}+ \sum_{n=1}^{\infty} \l\| \int\!\!\!\cdots\!\!\!\int_{\Delta_n(T_{2r}, t)} \nabla b\left(t_{n}, X_{t_{n}}^x\right)\cdots\nabla b\left(t_{1}, X_{t_{1}}^x\right) \ \d t_1 \cdots \d t_n \r\|_{\widetilde{L}^{pr}_xL^r_\om} \\ 
		&+ \sum_{n=1}^{\infty} \l\| \int\!\!\!\cdots\!\!\!\int_{\Delta_n(T_{2r}, t)} \nabla b\left(t_{n}, X_{t_{n}}^x\right)\cdots\nabla b\left(t_{1}, X_{t_{1}}^x\right) \ \d t_1 \cdots \d t_n \r\|_{{\widetilde{L}^{2pr}_xL^{2r}_\om}}  \|\nabla X_{T_{2r}}^x-\mathrm{I}\|_{\widetilde{L}^{2pr}_xL^{2r}_\om}\\
		\leq & C t^{\gamma/2r}. 
	\end{align*}
	By induction, we can show that the estimate in \eqref{eq-DX-Lr} holds for \(t \in [kT_{2r} \wedge T, (k+1)T_{2r} \wedge T]\). Thus, \eqref{eq-DX-Lr} holds for all \(t \in [0, T]\).
    
    Using \eqref{Eq-DXexp} one sees that \eqref{Eq-DX} can be proved in the same way as \eqref{Eq-nablaX}.  

    \medskip
    
    For \eqref{Eq-DX-Holder}.  Assume $0\leq s<s'\leq t\leq  T$. Combing \eqref{Eq-DXexp2} and \eqref{Eq-DX}, we obtain 
	\begin{align*}
		&\| D_{s} X_t^x -D_{s^{\prime}} X_t^x \|_{\widetilde L^{pr}_xL^r_\om} \\
		\leq&  \l\| D_{s} X_{s'}^x-{\rm I} \r\|_{\widetilde L^{pr}_xL^r_\om}+  \l\| D_{s} X_{s'}^x-{\rm I} \r\|_{\widetilde L^{2pr}_xL^{2r}_\om} \|D_{s'}X_t^x-{\rm I}\|_{\widetilde L^{2pr}_xL^{2r}_\om}\\
		\leq& C (s'-s)^{\gamma/4r}, 
	\end{align*}
	for a.e. $s,s'\in [0,T]$ with $0\leq s<s'\leq t\leq T$. 
	So, we complete our proof for the first case. 
	\
	\\
	{\bf Case (b).} 
	Let $q\in (2,q_1)$ and $\gamma=\frac{1}{q}-\frac{1}{q_1} \in (0,\frac12-\frac{1}{q_1})$. By \eqref{Eq-Keyest2-1} and the argument in the previous case, one can see that for each positive even integer $r$, 
	\const{\CXitetwo}
	\begin{align*}
		&\l\|\int\!\!\!\cdots\!\!\!\int_{\Delta_n(S_0,S_1)} \nabla b\left(t_{n}, X_{t_{n}}^x\right)\cdots\nabla b\left(t_{1}, X_{t_{1}}^x\right) \  \d t_{1} \cdots \d t_{n} \r\|_{\widetilde L^{pr}_xL^r_\om}\\
		\leq& (rC_{\CXitetwo})^{n+1/r} \|b\|_{\mL^{p_1}_{q_1}(S_0, S_1)}^{n} (S_1-S_0)^{\gamma/r}, 
	\end{align*}
	where $C_{\CXitetwo}>1$ only depends on $d, p_1, q_1, T,  a(m), r, p$ and $\gamma$.  Since $b\in {\mL^{p_1}_{q_1}(T)}$, for each even integer $r$ there is a positive  constant $T_r>0$ depending on $r, C_{\CXitetwo}$ and $\ell(\delta)$ such that for any $S_0, S_1\in [0,T]$ with $0\leq S_1-S_0\leq T_r$ 
	$$
	\|b\|_{\mL^{p_1}_{q_1}(S_0, S_1)} \leq (2rC_{\CXitetwo})^{-1}. 
	$$ 
	Thus, 
	\begin{align*}
		&\l\|\int\!\!\!\cdots\!\!\!\int_{\Delta_n(S_0,S_1)} \nabla b\left(t_{n}, X_{t_{n}}^x\right)\cdots\nabla b\left(t_{1}, X_{t_{1}}^x\right) \   \d t_{1} \cdots \d t_{n} \r\|_{\widetilde L^{pr}_xL^r_\om} \\
		\leq& (rC_{\CXitetwo})^{1/r} 2^{-n} (S_1-S_0)^{\gamma/r}. 
	\end{align*}
	Our desired estimates then can be obtained by the above estimate and the same argument  as in the previous case. \end{proof}

	\section{Proof of the main result}\label{Sec-Proof}
	The following lemma is a consequence of Theorem 1.1 in \cite{rockner2022weak}. 
	\begin{lemma}\label{Le-Weak}
	Let $d\geq 3$. Assume that $b\in C([0,T]; L^d)$ or $b\in \mL^{p_1}_{q_1}(T)$ with $p_1,q_1\in (2,\infty)$ and $d/p_1+2/q_1=1$. Then there is a unique weak solution to \eqref{Eq-SDE} such that for any $p,q\in (1,\infty)$ with $d/p+2/q<2$, the Krylov type estimate \eqref{Eq-Krylov} is valid. 
	\end{lemma}
	
	Now we are in the position to prove our main result.  
	\begin{proof}[Proof of Theorem \ref{Th-Main}] 
	{\bf Case (a): $b\in C([0,T]; L^d)$.} 
	\
	\\
	Recalling that $\rho\in C_c^\infty(\R^d)$  satisfying $\rho\geq 0$ and $\int \rho=1$, and let $b_k=b*_x\rho_k$. 
	Since $b\in C([0,T]; L^d)$, by Propostion \ref{Prop-CLd} we have 
	\begin{equation}\label{eq-K-bk}
	\begin{aligned}
		\| b_k-b_k*_x\rho_m \|_{\mL^d_\infty(T)} =& \| (b-b*_x\rho_m)*_x\rho_k \|_{\mL^d_\infty(T)}\\
		\leq& \| b-b*_x\rho_m \|_{\mL^d_\infty(T)}=:a(m)\to 0\quad (m\to\infty). 
	\end{aligned}
	\end{equation}
	It is well-known that for each $k$ there is a unique continuous random field $X(k): \Delta_2(T)\times \R^d\times \Om\to \R^d $ such that 
	\begin{equation}\label{Eq-mSDE}
	X^x_{s,t}(k)=x+\int_s^t b_k(r, X_{s,r}^x(k)) \d r+ W_t-W_s, \ \mbox{ for all }\  0\leq s\leq t\leq T,\  x\in \R^d. 
	\end{equation}
	Given $\beta\in (0, 1/2)$, let 
	$$
	p\in (1,d), \ q\in (1,\infty) \ \mbox{ satisfying } \ \frac{d}{p}+\frac{2}{q}\in (1,2-2\beta). 
	$$
	By estimate \eqref{eq-K-bk} and Remark \ref{Rek-local-W2}, for any $s\leq t_1\leq t_2\leq T$ and $f\in \tL^p_q(T)$, there is a unique function $u_k$ in $\tH^{2,p}_q(T)$ solving 
	$$
	\p_t u_k + \frac{1}{2}\Delta u_k + b_k \cdot \nabla u_k + f=0 \mbox{ in } (s, t_2)\times \R^d, \quad u_k(t_2)=0
	$$
	and a constant $C$, which does not depends on $k$, such that 
	\begin{equation}\label{eq-uk-W2}
	\|\p_t u_k\|_{\tL^{p}_q(t_1,t_2)} + \|u_k\|_{\tH^{2,p}_q(t_1,t_2)} \leq C \|f\|_{\tH^{p}_q(t_1,t_2)}.  
	\end{equation}
	By the generalized It\^o formula (cf.  \cite{rockner2022weak}), 
	$$
	-u_k(t_1, X_{s,t_1}^x(k))= - \int_{t_1}^{t_2} f(t, X_{s,t}^x(k)) \ \d t + \int_{t_1}^{t_2}\nabla u_k(t, X_{s,t}^x(k)) \cdot \d W_t. 
	$$
	
	Taking $\a=0$ and $\theta=1-\frac{1}{q}-\beta>\frac{d}{2p}$ in \eqref{Eq-PMorrey}, using Morrey's inequality and \eqref{eq-uk-W2} we get 
	\begin{equation}\label{eq-krylov1}
	\begin{aligned}
		&\bE \l(\int_{t_1}^{t_2} f(t, X_{s,t}^x(k)) \d t \Big| \sF_{t_1}\r)=\bE \l( u_k(t_1, X_{s,t_1}^x(k))\Big |\sF_{t_1}\r)\\
		\leq& \|u_k(t_1)\|_{\infty} \leq C \sup_{z\in \R^d}\|(u_k\chi_1^z)(t_1)\|_{H^{2\theta, p}}\\\overset{\eqref{Eq-PMorrey}}{\leq} &  C  |t_2-t_1|^{\beta}  \l( \|\p_t u_k\|_{\tL^{p}_q(t_1,t_2)} + \|u_k\|_{\tH^{2,p}_q(t_1,t_2)} \r)\\\overset{\eqref{eq-uk-W2}}{\leq} & 
		C |t_2-t_1|^{\beta} \|f\|_{\tL^{p}_{q}(t_1, t_2)},   
	\end{aligned}
	\end{equation}
	where $C$ only depends on $d, p, q, \beta, T, \|b\|_{\mL^d_\infty(T)}$ and $\{a(m)\}$. Once with \eqref{eq-krylov1} in hand, it is standard to show that 
	\begin{equation}\label{eq-krylov2} 
	\bE \l| \int_{t_1}^{t_2} f(t, X_{s,t}^x(k)) \d t \r|^r \leq C_r |t_2-t_1|^{\beta r} \|f\|_{\tL^{p}_{q}(t_1, t_2)}^r, \quad \forall \ f\in \tL^p_q(T), \ r>0, 
	\end{equation}
	(cf. \cite{zhang2018singular}). Therefore, by noting that $p<d$, $q<\infty$ and $\beta<1/2$, we get 
	\begin{equation}\label{Eq-Xt1Xt2}
	\begin{aligned}
		\bE \l|X_{s, t_1}^x(k)-X_{s, t_2}^x(k)\r|^r\leq& C \bE \l(\int_{t_1}^{t_2} |b_k|(t, X_{s,t}^x) \d t \r)^r +C \bE | W_{t_2}-W_{t_1}|^r\\
		\leq& C |t_2-t_1|^{\beta r} \l( 1+ \|b\|_{\mL^{d}_{\infty}(T)}^r\r), \quad \forall r> 0. 
	\end{aligned}
	\end{equation}
    Consequently, 
    \begin{equation}\label{eq-X-Lr}
	\sup_{\substack{x\in \R^d;\\ 0\leq s\leq t\leq T}} \bE \l|X_{s, t}^x(k)-x\r|^r\leq  C, \quad \forall r>0.  
    \end{equation}
    On the other hand, \eqref{Eq-nablaX} and Morrey's inequality implies  that for each $r>d$, 
    \begin{align*}
    \sup_{\substack{z\in \R^d; \\0\leq s\leq t\leq T}}\bE \|X_{s, t}^x(k)-x\|^{r}_{\dot{C}_x^{1-\frac{d}{r}}(B_1(z))} \leq C \sup_{\substack{z\in \R^d; \\0\leq s\leq t\leq T}}\bE  \|\nabla X_{s, t}^x(k)-\mathrm{I}\|_{L^{r}_x(B_1(z))}^r\leq C. 
    \end{align*}
    Therefore, for any $0\leq s\leq t\leq T,  x, y\in \R^d$,  
	\begin{equation}\label{Eq-X(x)X(y)}
	\begin{aligned}
		\bE |X_{s, t}^x(k)-X_{s, t}^y(k))|^r \leq C |x-y|^{r-d}+C|x-y|^r, \quad \forall r>d, 
	\end{aligned}
	\end{equation}
	where $C$ only depends on $d, p, q, r, \|b\|_{\mL^d_\infty(T)}$ and $\{a(m)\}$. 
	
	Assume $0\leq s_1\leq s_2\leq t$. By the Markov property and the independence of $X_{s_1, s_2}^x(k)$ and $X_{s_2, t}^y(k)$, for each $r>d$ we obtain 
	\begin{equation}\label{eq-Xs1Xs2}
	\begin{aligned}
		&\bE  |X_{{s_1}, t}^x(k)-X_{s_2, t}^x(k)|^r\\
		\leq& C_r \bE \l|\int_{s_1}^{s_2} b_k \l(s, X_{s_1, s}^x(k)\r) \d s\r|^r + C_r \bE \l| \int_{s_2}^t \l[ b_k\l(s, X_{s_1, s}^x(k)\r)- b_k \l(X_{s_2, s}^x(k)\r) \r] \d s \r|^r\\
		\overset{\eqref{eq-krylov2}}{\leq} & C |s_1-s_2|^{\beta r} + C \bE  \l| \int_{s_2}^t \l[b_k(s, X_{s_2, s}^{X_{s_1, s_2}^x(k)}(k))- b_k(s, X_{s_2, s}^x(k))\r] \d s\r|^r\\ 
		\leq & C |s_1-s_2|^{\beta r} + C \bE  \l| X_{s_2, t}^{X_{s_1, s_2}^x(k)}(k)- X_{s_2, t}^x(k)\r|^r\\
		=& C |s_1-s_2|^{\beta r} + C \bE  \l[ \bE \l|X_{s_2, t}^y(k)-X_{s_2, t}^x(k)\r|^r \Big|_{y=X_{s_1, s_2}^x(k)}\r]\\
		\overset{\eqref{Eq-X(x)X(y)}}{\leq}& C |s_1-s_2|^{\beta r} + C \bE  \l|X_{s_1, s_2}^x(k)-x\r|^{r-d}+C\bE  \l|X_{s_1, s_2}^x(k)-x\r|^{r} \\
		\overset{\eqref{Eq-Xt1Xt2}}{\leq} & C |s_1-s_2|^{\beta (r-d)} . 
	\end{aligned}
	\end{equation}
	Combing \eqref{Eq-Xt1Xt2}, \eqref{Eq-X(x)X(y)} and \eqref{eq-Xs1Xs2}, we obtain that for all $(s_i, t_i)\in \Delta_2(T)$, $i=1,2$, 
 .  \begin{equation}\label{eq-holder-Xn}
    \begin{aligned}
        &\bE |X_{s_1, t_1}^x(k)-X_{s_2, t_2}^y(k))|^{r}\\
        \leq& C \l(|t_1-t_2|^{\beta r}+ |x-y|^{r-d} + |x-y|^{r} + |s_1-s_2|^{\beta (r-d)}\r), 
    \end{aligned}
    \end{equation}
	where $r>d$ and $C$ only depends on $d, r, T$ and $b$. On the other hand, noting that 
	$$
	\sup_k \| b_k  \|_{\mL^d_\infty(T)} \leq \| b \|_{\mL^d_\infty(T)}, \quad \| b_k-b_k*_x\rho_m \|_{\mL^d_\infty(T)} \leq a(m)\to 0\ (m\to \infty),  
	$$
	by Lemma \ref{Le-Comp} and Proposition \ref{Prop-Key} one can see that for any fixed $(s, t)\in \Delta_2(T)$ and $R>0$, 
    $$
        \l\{ B_R\times \Om \ni (x, \om)\mapsto X^x_{s, t}(k)(\om) \in \R^d \r\}_{k\in \mN_+}
    $$
    is relatively compact in $L^2(B_R\times \Om)$. The standard diagonal argument yields that there is a subsequence (still denoted by $X^x_{s,t}(k)$) and a countable dense subset $\cD$ of $\R^d$ such that 
	$$
	X_{s,t}^x(k)\xrightarrow[k\to\infty]{L^2(\Om) \mbox{ {\tiny and} } a.s.} X_{s, t}^x, \ \forall (s,t) \in \mQ^2\cap \Delta_2(T) \mbox{ and }x\in \cD. 
	$$
	By \eqref{eq-X-Lr}, we also have 
	$$
	X_{s,t}^x(k)\xrightarrow[k\to\infty]{L^r(\Om)} X_{s, t}^x, \ \forall r\geq 1,\, \forall (s,t) \in \mQ^2\cap \Delta_2(T) \mbox{ and }x\in \cD. 
	$$ 
	Fatou's lemma and \eqref{eq-holder-Xn} yield that for all $(s_i,t_i) \in \mQ^2\cap \Delta_2(T)$, $i=1,2$, and $x\in \cD$, 
    \begin{equation}\label{eq-holder}
    \begin{aligned}
        &\bE \l|X_{s_1, t_1}^{x_1}-X_{s_2, t_2}^{x_2}\r|^{r} \\
        \leq& C \l(|x_1-x_2|^{r-d} + |x_1-x_2|^{r}+|s_1-s_2|^{\beta (r-d)}+ |t_1-t_2|^{\beta r}\r), \quad \forall r>d. 
    \end{aligned}	
    \end{equation}
    
    Therefore, $X_{s,t}^x$ can be extended to a continuous random field on $\Delta_2(T)\times \R^d$ satisfying \eqref{eq-holder} due to the Kolmogorov-Chentsov theorem, and up to a subsequence (still denoted by $X_{s,t}^x(k)$), 
	\begin{equation}\label{eq-XktoX}
	X_{s,t}^x(k,\om) \overset{k\to\infty}{\longrightarrow}  X_{s,t}^x(\om), 
	\end{equation}
	for all $(s,t)\in \mQ^2\cap \Delta_2(T)$, $x\in \cD$ and $\bP$-a.s. $\om\in \Om$. Then, since by \eqref{eq-holder-Xn} for $\bP$-a.s. $\om\in \Om$, $X^x_{s,t}(k, \om)$, $k\in \mN_+$ are equicontinuous as functions of $(s,t,x)$, \eqref{eq-XktoX} holds for all $(s,t) \in \Delta_2(T), x\in \R^d$ and $\om\in \Om_0\in \sF$ with $\bP(\Om_0)=1$. Moreover, by Proposition \ref{Prop-Key} and the definition of \((X_{s,t}^x)\), it satisfies \eqref{Eq-gradient}. 
    
    Taking limits on both sides of \eqref{eq-krylov2}, we get  
	\begin{equation}\label{eq-krylov3}
	\bE \l|\int_{t_1}^{t_2} f(t, X_{s,t}^x) \d t \r|^r\leq C |t_2-t_1|^{\beta r} \|f\|_{\tL^{p}_{q}(t_1, t_2)}^r.  
	\end{equation} 
	Thus, for each $x\in \R^d$ and $K\in \mN_+$, 
	\begin{align*}
		&\bE \sup_{t\in [s,T]}\l|\int_{s}^t b(\tau, X_{s, \tau}^x) \d \tau-\int_{s}^t b_{k}(\tau, X_{s, \tau}^x(k)) \d \tau \r|\\
		\leq& \bE \int_{s}^T \l| b-b_{K}\r|(\tau, X_{s, \tau}^x) \d \tau+  \bE \int_{s}^T \l| b_{K}-b_{k}\r|(\tau, X_{s, \tau}^x(k)) \d \tau \\ 
		&+\bE\sup_{t\in [s,T]}\l|\int_{s}^t b_{K}(\tau, X_{s, \tau}^x) \d \tau - \int_{s}^t b_{K}(\tau, X_{s, \tau}^x(k)) \d \tau \r|. 
	\end{align*}
	By our assumption on $b$, it holds that $b-b_k\to 0$ in $\mL^d_\infty(T)$. So, the first and second terms on the right-hand side of the above inequality converge to $0$ as $k$ goes to infinity, due to the fact that $X$ and $X(k)$ satisfy the Krylov type estimates \eqref{eq-krylov2} and \eqref{eq-krylov3}. On the other hand, by \eqref{eq-XktoX} and Lebesgues dominated convergence theorem, the third term on the right side of the above inequality also converges to $0$ as $k$ goes to infinity. So, 
	$$
	\bE\sup_{t\in [s,T]}\l|\int_{s}^t b(\tau, X_{s, \tau}^x) \d \tau-\int_{s}^t b_{k}(\tau, X_{s, \tau}^x(k)) \d \tau \r|\to 0, 
	$$
	which together with \eqref{eq-XktoX} implies 
	$$
	X_{s,t}^x-x-\int_s^t b(\tau, X_{s,\tau}^x) \d \tau= \lim_{k\to\infty}\l(X_{s,t}^x(k)-x-\int_s^t b(\tau, X_{s,\tau}^x(k)) \d \tau\r)=W_t,
	$$
	i.e. the limit point $X^x_{s,\cdot}$ is a strong solution to \eqref{Eq-SDE}. Hence, we obtain the strong existence of solutions to \eqref{Eq-SDE}. Moreover, by the proof of  Theorem of 1.1 in \cite{rockner2022weak}, we can also see that $X^x_{s,\cdot}$ also satisfies \eqref{Eq-Krylov} for any $p', q'\in (1,\infty)$ satisfying $d/p'+2/q'<2$. 
	
	Following \cite{cherny2002uniqueness}, we next show that the limit point of $X_{s,t}^x(k)$ is the unique, and is also the unique strong solution to \eqref{Eq-SDE} satisfying \eqref{Eq-Krylov}. Without loss of generality we may assume $s=0$. Suppose $X$ is a limit point of $X^x_{\cdot}(k)$, which is a strong solution of \eqref{Eq-SDE} with $s=0$ on $(\Om, \sF, \{\sF_t\}_{t\in [0,T]}, \bP)$. Then there exists a measurable map $\cT: C([0, T];\R^d)\to C([0, T];\R^d)$ such that $X_{\cdot}(\om)=\cT(W(\om))$ for $\bP$-a.s. $\om$. Let $\{\bQ_\om\}_{\om\in \Om}$ be the regular conditional expectation of $X$ with respect $ \sF^W_{T} := \sigma\{ W_t: t\in [0,T]\}$. Then $\bQ_\om=\delta_{\cT(W(\om))}$ for $\bP$-a.s. $\om$. Now, let $Y$ be an another strong solution to \eqref{Eq-SDE} with $s=0$ on $(\Om, \sF, \{\sF_t\}_{t\in [0,T]}, \bP)$ satisfying \eqref{Eq-Krylov}. Thanks to Lemma \ref{Le-Weak}, we have $\mathrm{law}(X)=\mathrm{law}(Y)$, together with the fact that 
	$$
	W_t= X_{t}-x-\int_0^t b(r, X_{r}) \d r= Y_{t}-x-\int_0^t b(r, Y_{r}) \d r, 
	$$ 
	we obtain $\mathrm{law}(X, W)=\mathrm{law}(Y,W)$. This implies $\bQ'_\om$, the regular conditional expectation of $Y$ with respect to $\sF^W_{T}$,  equals to $\bQ_\om$ for $\bP$-a.s. $\om$, i.e. $\bQ'_\om=\delta_{\cT(W(\om))}$. Thus, $Y(\om)=\cT(W(\om))=X(\om)$ for $\bP-a.s.\ \om$. 
    
    Moreover, the estimates in \eqref{Eq-gradient} and \eqref{Eq-Holder-X} follow from Proposition \ref{Prop-Key} and \eqref{eq-holder}, respectively.  
	
	{\bf Case (b): $b\in \mL^{p_1}_{q_1}(T)$.} Given $\beta\in (0,1/2)$. In this case, we take $p\in (1,p_1)$ and $q\in (1, q_1)$ such that $d/p+2/q\in (1, 2-2\beta)$. Define the maximal function of $b(t,\cdot)$: 
	$$
	\cM b(t,x): =\sup_{r>0}\fint_{B_r(x)} \l|b(t, y)\r| \d y
	$$
	Define $b_k= (b\1_{|b|\leq k})*_x\rho_k\in L^\infty([0,T]; C^2_b)$. 
	Noting that $b_k\leq |b|*_x\rho_k\leq C\cM b$ (cf. \cite[Corollary 2.8]{duoandikoetxea2001fourier}), we have $K'_{b_k}(m) \leq K'_{C \cM b}(m)$. By the basic fact that 
	$$\|\cM b\|_{\mL^{p_1}_{q_1}(T)}\asymp \|b\|_{\mL^{p_1}_{q_1}(T)}<\infty, 
	$$
	(cf. \cite[Theorem 2.5]{duoandikoetxea2001fourier}), we obtain 
	$$
	\sup_{k} K'_{b_k}(m) \leq K'_{C \cM b}(m)=: a(m) \to 0, \mbox{ as } m\to \infty. 
	$$
	Also we have 
	$$
	\sup_{k} \om_{b_k}(\delta)\leq \om_b(\delta)=:\ell(\delta)\to 0, \mbox{ as } \delta\to 0. 
	$$
	Then our desired results in the second case can be obtained by the same procedure as for the previous case. 
	\end{proof}
	\begin{remark}
	By Remark \ref{Rek-Global}, if we further assume $b\in \mL^p_{q}(T)$ with $p\in (\frac{d}{d-1},d)$,  $q=\gamma^{-1}$ in case (a) and $p\in (\frac{p_1}{p_1-1}, p_1)$,  $q=(1/q_1+\gamma)^{-1}$ in case (b), then 
	\begin{equation}\label{Eq-gradient-1}
	\sup_{\substack{0\leq s\leq t\leq T}} \int_{\R^d}\l(\bE |\nabla X_{s,t}^x-\mathrm{I} |^r \r)^p\d x <\infty, \ \mbox{ for any } r\in [2,\infty). 
	\end{equation}
	\end{remark}

	\section{Application}\label{Sec-Application}
	
	In \cite{zhang2010astochastic}, Zhang studied the backward Navier-Stokes equation \eqref{Eq-BNSE} 
	through considering the stochastic system \eqref{Eq-SS}.  As in \cite{constantin2008stochastic}, it was also  shown in \cite{zhang2010astochastic} that the existence of smooth solutions for \eqref{Eq-BNSE} and \eqref{Eq-SS} are equivalent (see Theorem 2.3 therein). Therefore, it is quite interesting to find a regularity criterion for solutions of \eqref{Eq-SS}. Below we give one such conditional regularity result, which is similar to the Serrin criterion for the Navier-Stokes equations.  
	
	\begin{theorem}\label{Th-Main2}
	Let $d\geq 3$, $T>0$, $p_1, q_1\in (2,\infty)$, $q>d$ and $k,l\in \mN$. Assume $u\in C([-T,0];L^d)\cap \mL^2_\infty(-T,0)$ or $u\in \mL_{q_1}^{p_1}(-T,0)\cap\mL^2_\infty(-T,0)$ with $d/p_1+2/q_1\leq 1$ and  $\varphi\in H^{k, q}$. Suppose that $(u, X)$ is a solution to the stochastic system \eqref{Eq-SS} and $X$ satisfying the Krylov type estimate \eqref{Eq-Krylov}, then $u\in \mH^{k,q}_\infty(-T,0)$ and for any $l\leq k/2$, $\p_t^{l} u\in \mL^q_\infty(-T,0)$. Consequently, if $\varphi\in C_c^\infty(\R^d)$, then $u\in C^\infty_b([-T, 0]\times\R^d)$ and it satisfies \eqref{Eq-BNSE}. 
	\end{theorem}
	
	\begin{proof}
	{\em Step 1.} Assume that $\varphi\in L^q$,  $u\in C([-T,0]; L^d)\cap \mL^2_\infty(-T,0)$, or $u\in \mL^{p_1}_{q_1}(-T,0)\cap \mL^2_\infty(-T,0)$ with $p_1, q_1\in (2,\infty)$ and $d/p_1+2/q_1\leq 1$. We claim that
	\begin{equation}\label{eq-u-Lq}
	\sup_{t\in [-T,0]} \|u(t)\|_q<\infty.  
	\end{equation}
	Below we only give the proof of \eqref{eq-u-Lq} for the case where  $u\in \mL^{p_1}_{q_1}(-T,0)\cap\mL^2_\infty(-T,0)$ with $p_1, q_1\in (2,\infty)$ and $d/p_1+2/q_1= 1$, since the first case is simpler. 
	
	Fix $a\in (d, q)$ and define 
	$$
	v(t,x):= \mathbf{E}\left[( \nabla^{\top}  X_{t,0}^x -{\rm I})\varphi\left( X_{t, 0}^x \right)\right]. 
	$$
	Noting that $\frac{1}{a}-\frac{1}{q}\in (0, \frac{1}{d})$, one can always choose $r\geq 2$ and $p\in (\frac{p_1}{p_1-1}, p_1)$ such that $\frac{1}{pr}=\frac{1}{a}-\frac{1}{q}$. Thus, for each $t\in [-T, 0]$ 
	\begin{equation}\label{eq-v-La}
		\begin{aligned}
			\|v(t)\|_{a} \leq & \Big\| \l\|\nabla^{\top}  X_{t,0}^x-{\rm I}\r\|_{L^{r}_\om}  \l\|\varphi(X_{t,0}^x)\r\|_{L^{r'}_\om}  \Big\|_{L^a_x} \\
			\leq &  C   \l\|\nabla X_{t,0}^x-{\rm I}\r\|_{L^{pr}_xL^{r}_\om}  \l\|\varphi(X_{t,0}^x)\r\|_{L^{q}_x L^{r'}_\om}\\ 
			\overset{\eqref{Eq-gradient-1}}{\leq} & C \| \bE |\varphi|^{r'}(X_{t,0}^x)\|_{L^{q/r'}_x}^{1/r'}\leq C \|\varphi\|_{q}. \end{aligned}
	\end{equation}
	Here $r'=r/(r-1)$, and we use the fact that 
	\begin{equation}\label{eq-lp-persistent}
	\|\bE f(X^x_{t, 0})\|_{L^q_x} \leq \|f\|_q, \quad \forall q\in [1,\infty],  
	\end{equation}
	due to the fact that $u$ is divergence free (cf. \cite[Lemma 3.2]{zhang2020stochastic}). 
	Recall that $\mathrm{P}$ is the Leray projection 
	$$
	(\mathrm{P}F)_i=F_i-\nabla (\Delta)^{-1}\div F = F_i-\sum_{j=1}^d R_iR_j F_j, 
	$$ 
	where $R_i$ is the Riesz transformation. The $L^q$ boundedness of $R_i$ implies that $\mathrm{P}$ is a bounded map on $L^q(\R^d; \R^d)$ with $q\in (1,\infty)$. By \eqref{eq-lp-persistent}, we have 
	\begin{align*}
		\|u\|_{\tL^a_\infty(-T,0)}=&\sup_{z\in \R^d} \|u\chi_1^z\|_{\mL^a_\infty(-T,0)}\\
		\leq & C\l (\sup_{t\in [-T, 0]}  \| {\rm P} v(t)\|_{a} + \sup_{t\in [-T, 0]}  \|{\rm P} \bE \varphi (X^\cdot_{t,0}) \|_q \r)\\
		\leq& C  \sup_{t\in [-T, 0]} \|v(t)\|_a + C  \|\varphi\|_{q}\overset{\eqref{eq-v-La}}{<}\infty. 
	\end{align*}
	Noting that $a>d$, combining the above estimate and Theorem 1.1 of \cite{xia2020lqlp} we get 
	\begin{equation}\label{eq-X'-Lr}
	\sup_{x\in \R^d}\bE \sup_{t\leq s\leq 0} |\nabla X_{t, s}^x |^r<\infty, \quad \forall r\geq 1. 
	\end{equation}
	Therefore, for each $t\in [-T,0]$, 
	\begin{equation}
	\begin{aligned}
		\|u(t)\|_{q} \leq & \Big\| \l\|\nabla^{\top}  X_{t,0}^x\r\|_{L^{r}_\om}  \l\|\varphi(X_{t,0}^x)\r\|_{L^{r'}_\om}  \Big\|_{L^q_x} \leq  C   \l\|\nabla X_{t,0}^x\r\|_{L^{\infty}_xL^{r}_\om}  \l\|\varphi(X_{t,0}^x)\r\|_{L^{q}_x L^{r'}_\om}\\ 
		\overset{\eqref{eq-X'-Lr}}{\leq} & C \| \bE |\varphi|^{r'}(X_{t,0}^x)\|_{L^{q/r'}_x}^{1/r'}\leq C \|\varphi\|_{q}<\infty.    
	\end{aligned}
	\end{equation} 
	So, we complete our proof for \eqref{eq-u-Lq}. 
	
	{\em Step 2.} Now assume that $\varphi\in H^{1,q}$. By \cite[Lemma 7.2]{zhang2016stochastic}, 
	\begin{equation}\label{eq-pi-u}
	\begin{aligned}
		\p_i u(t)= \p_i \mathrm{P} \bE \left[\nabla^{\top} X_{t,0}^x \varphi(X_{t, 0}^x) \right] 
		=& \mathrm{P}\bE \left[\nabla^{\top} X_{t, 0}^x [\nabla \varphi (X_{t, 0}^x)-\nabla^{\top} \varphi (X_{t, 0}^x)]\p_i X_{t, 0}^x\right].
	\end{aligned}
	\end{equation}
	Using  H\"older's inequality and \eqref{eq-X'-Lr}, we get 
	\begin{equation*}
		\begin{aligned}
			\|\p_i u(t)\|_q\leq &C \l\|\bE \left[\nabla^{\top} X_{t, 0}^x [\nabla \varphi (X_{t, 0}^x)-\nabla^{\top} \varphi (X_{t, 0}^x)]\p_i X_{t, 0}^x\right] \r\|_{L^q_x}\\
			\leq & C  \Big\| \|\nabla X^x_t \|_{L^{2r}_\om} \|\nabla \varphi (X_{t, 0}^x)\|_{L^{r'}_\om}\Big\|_{L^q_x}\\
			\leq & C \|\nabla X^x_t \|_{L^{\infty}_xL^{2r}_\om} \|\nabla \varphi (X_{t, 0}^x)\|_{L^{q}_xL^{r'}_\om}\overset{\eqref{eq-X'-Lr}}{\leq} C \|\nabla \varphi \|_{q}<\infty. 
		\end{aligned}
	\end{equation*}
	Hence, 
	\begin{equation}\label{eq-u-H1}
	\|u\|_{\mH^{1,q}_\infty(-T,0)} \leq C \| \varphi \|_{H^{1,q}}<\infty.  
	\end{equation} 
	
	{\em Step 3.} Assume that $\varphi\in H^{2,q}$. Following \cite{xia2020lqlp}, below we use a Zvonkin type change of variables to convert the first equation in \eqref{Eq-SS} to a new SDE. 
	Let $t\in [-T, 0]$, $\lambda\geq 0$ and $a\in (1,\infty)$ such that $d/q+2/a<1$. Since $u\in \mH^{1,q}_\infty(-T,0)$, there is a unique function $U$ in $\mH^{3, q}_a(-T,0)$ satisfying 
	$$
	\p_s U + \l(\frac{\Delta}{2}-\lambda\r)U+ u\cdot \nabla U+u=0 \ \mbox{ in }(-T,0)\times\R^d, \quad U(0)=0. 
	$$
	Moreover, 
	\begin{equation}\label{eq-U-W3}
	\lambda \|U\|_{\mH^{1,q}_a(-T,0)}+\|\p_t U\|_{\mH^{1,q}_a(-T,0)} +  \|\nabla^2 U\|_{\mH^{1,q}_a(-T,0)} \leq C \|u\|_{\mH^{1,q}_a(-T,0)}<\infty
	\end{equation}
	(cf. \cite{xia2020lqlp}). Since $d/q+2/a<1$, using \eqref{eq-U-W3}, \eqref{Eq-PMorrey} and an interpolation inequality one can choose $\lambda$ large enough so that 
	\begin{equation}\label{eq-U'}
	\sum_{k=0}^2 \|\nabla^k U\|_{\infty}\leq 1/2.  
	\end{equation}
	Define 
	$$
	\Phi(s,x):= x+ U(s, x). 
	$$
	By \eqref{eq-U'}, $\Phi(s, \cdot)$ is a $C^2$-diffeomorphism and 
	\begin{equation}\label{eq-Phi-Phi-1}
	\| \nabla \Phi\|_\infty,\| \nabla^2 \Phi\|_\infty, \| \nabla \Phi^{-1}\|_\infty, \| \nabla^2 \Phi^{-1}\|_\infty \leq C.
	\end{equation}
	Set 
	$$
	Y_{t,s}^{y,k}:= \Phi^k(s, X_{t, s}^{\Phi^{-1}(t,y)}),\quad  \sigma^k_{k'}(s,y)= \p_{k'} \Phi^{k}(s,\Phi^{-1}(s,y)), \quad b^k(s, y)=\lambda U^k(s, \Phi^{-1}(s,y)). 
	$$
	Then, 
	$$
	Y_{t,s}^y=y + \int_t^s b(\tau , Y_{t,\tau }^y) \d \tau + \int_t^s \sigma(\tau , Y_{t,\tau }^y)  \d \widetilde W_\tau , 
	$$
	where $\widetilde W_\tau := W_\tau -W_t$ is a standard Brownian motion on $[t, 0]$. By \eqref{eq-U-W3}-\eqref{eq-Phi-Phi-1}, \eqref{Eq-PMorrey} and the definitions of $\si$ and $b$, one sees that 
	\begin{equation}\label{eq-sigma}
	\si(s,y)-{\rm I}=\nabla U(s, \Phi^{-1}(s,y))\in \mH^{2,q}_a(-T,0)\cap C([-T, 0]; C_b^1(\R^d))
	\end{equation}
	and 
	\begin{equation}\label{eq-b}
	b\in C([-T,0];C^2_b(\R^d)). 
	\end{equation}
	By the proof for \cite[Theorem 1.1]{xia2020lqlp}, $\p_i Y^y_{t,s}$ satisfies 
	\begin{equation}\label{eq-Y'}
	\p_{i} Y^{y}_{t,s}=e_i+\int_t^s \nabla b (\tau , Y_{t,\tau }^y) \, \p_{i} Y_{t,\tau }^{y} \,\d \tau + \int_t^s \p_l \si_{k'} (\tau , Y_{t,\tau }^y) \, \p_{i} Y_{t,\tau }^{y,l}\, \d \widetilde W_\tau ^{k'}.  
	\end{equation}
	and 
	\begin{equation}\label{eq-Y'-Lr}
	\sup_{x\in \R^d}\bE \sup_{t\leq s\leq 0} |\nabla Y_{t, s}^x |^r<\infty, \quad \forall r\geq 1. 
	\end{equation}
	Noting that \( \nabla b \in C([-T,0]; C^2_b(\mathbb{R}^d)) \) and \( \nabla \sigma \in \mathbb{H}^{1,q}_a(-T,0) \), one can again follow the main argument in the proof of \cite[Theorem 1.1]{xia2020lqlp} to see that 
	$$
	\sup_{x\in \R^d}\bE\sup_{s\in [t,0]} |\nabla^{2} Y_{t, s}^x |^r <\infty,  \quad \forall r\geq 1. 
	$$
	and
	\begin{equation*}
		\begin{aligned}
			\p_{ij} Y^{y}_{t,s}=& \int_t^s \p_{ll'}b(\tau , Y_{t,\tau }^y) \, \p_i Y_{t,\tau }^{y, l} \, \p_j Y_{t, \tau }^{y, l'} \,\d \tau  + \int_t^s \p_l b (\tau , Y_{t,\tau }^y) \, \p_{ij} Y_{t,\tau }^{y,l} \,\d \tau \\
			&+ \int_t^s \p_{ll'}\sigma_{k'}(\tau , Y_{t,\tau }^y) \, \p_i Y_{t,\tau }^{y, l} \, \p_j Y_{t, \tau }^{y, l'}\, \d \widetilde W_\tau ^{k'} + \int_t^s \p_l \si_{k'} (\tau , Y_{t,\tau }^y) \, \p_{ij} Y_{t,\tau }^{y,l}\, \d \widetilde W_\tau ^{k'}
		\end{aligned}
	\end{equation*}
	Recalling that $\Phi^{-1} \in C([-T, 0];C^2_b(\R^d))$, we see that  
	\begin{equation}\label{eq-D2X-Lr}
	\sup_{x\in \R^d}\bE\sup_{s\in [t,0]} |\nabla^{2} X_{t, s}^x |^r <\infty,  \quad \forall r\geq 1. 
	\end{equation}
	By \eqref{eq-pi-u}, we have 
	\begin{align*}
		\p_{ij} u(t)= \p_{ij} \mathrm{P} \bE \left[\nabla^{\top} X_t^x \varphi(X_t^x) \right]= \mathrm{P}  \bE \p_j \left[\nabla^{\top} X_t^x [\nabla \varphi (X_t^x)-\nabla^{\top} \varphi (X_t^x)]\p_i X_t^x\right].  
	\end{align*} 
	Using \eqref{eq-D2X-Lr} and following the same procedure as in the proof for \eqref{eq-u-H1}, one can verify that  
	$$
	\|u\|_{\mH^{2,q}_\infty(-T,0))}\leq C \| \varphi \|_{H^{2,q}}<\infty.
	$$
	Repeating the above process higher derivatives can be estimated similarly step by step. 
	
	{\em Step 4.} 
	Assume $\varphi \in H^{2,q}$ and set 
	$$
	w(t):= \bE \l[\nabla^\top X_{t,0}^x \, \varphi(X^x_{t,0}) \r]. 
	$$
	By \em Step 3 we can see that 
	$$
	\sup_{x\in \R^d}\bE \sup_{t\leq s\leq 0} |\nabla^3 X_{t, s}^x |^r<\infty \ \mbox{ and }\ 
	w\in \mH^{2,q}_\infty(-T,0). 
	$$
	Following the proof for \cite[Theorem 2.1]{zhang2010astochastic}, we see that $w$ satisfies 
	\begin{equation}\label{eq-w}
	\p_t w =- \tfrac{\Delta}{2} w- (\nabla w )u- (\nabla^\top u) w, \quad w(0)=\varphi.  
	\end{equation}
	Thus, $\p_t w\in \mL^q_\infty(-T, 0)$, which also implies $\p_t u=\p_t {\rm P} w={\rm P}\p_t w\in  \mL^q_\infty(-T, 0)$ due to the $L^q$ boundedness of ${\rm P}$. 
	
	If $\varphi \in H^{4, q}$, following the above discussion we see that $u, w\in \mH^{4,q}_\infty(-T, 0)$, which implies that the right side of \eqref{eq-w} is in $\mH^{2,q}_\infty(-T,0)$. Hence, $\p_t w\in \mH^{2,q}_\infty(-T,0)$ and $\p_t u\in \mH^{2,q}_\infty(-T,0)$. This means that $\p_t [\tfrac{\Delta}{2} w+ (\nabla w )u+(\nabla^\top u) w]\in \mL^q_\infty(-T,0)$, i.e. $\p_t^2 w \in \mL^q_\infty(-T,0)$. Repeating the same process one sees that $\p_t^k w \in \mL^q_\infty(-T, 0)$, provided that $\varphi\in H^{2k, q}$. So, we complete our proof.

	\end{proof}

	\section*{Acknowledgement}
	
	This work was carried out while the second author was employed at the University of Bielefeld, during which time he received sponsorship from the German Research Foundation (DFG) through the Collaborative Research Centre (CRC) ``Taming uncertainty and profiting from randomness and low regularity in analysis, stochastics and their applications,"-Project-ID 317210226-SFB 1283. The second named author is very grateful to Professor Nicolai Krylov and Xicheng Zhang who encouraged him to persist in studying this problem, and also to Professor Ka{\ss}mann for providing him with an excellent environment to work at Bielefeld University. 

	\newcommand{\etalchar}[1]{$^{#1}$}


\begin{thebibliography}{MPMBN{\etalchar{+}}13}
		
		\bibitem[BFGM19]{beck2019stochastic}
		Lisa Beck, Franco Flandoli, Massimiliano Gubinelli, and Mario Maurelli.
		\newblock Stochastic {ODE}s and stochastic linear {PDE}s with critical drift:
		regularity, duality and uniqueness.
		\newblock {\em Electronic Journal of Probability}, 24:1--72, 2019.
		
		\bibitem[BS04]{bally2004relative}
		Vlad Bally and Bruno Saussereau.
		\newblock A relative compactness criterion in {Wiener--Sobolev} spaces and
		application to semi-linear stochastic {PDE}s.
		\newblock {\em Journal of Functional Analysis}, 210(2):465--515, 2004.
		
		\bibitem[Che02]{cherny2002uniqueness}
		Aleksander~Semenovich Cherny.
		\newblock On the uniqueness in law and the pathwise uniqueness for stochastic
		differential equations.
		\newblock {\em Theory of Probability \& Its Applications}, 46(3):406--419,
		2002.
		
		\bibitem[CI08]{constantin2008stochastic}
		{P}eter {C}onstantin and {G}autam {I}yer.
		\newblock {A} stochastic {L}agrangian representation of the three-dimensional
		incompressible {{N}}avier-{S}tokes equations.
		\newblock {\em {C}ommunications on {P}ure and {A}pplied {M}athematics},
		61(3):330--345, 2008.
		
		\bibitem[Dav07]{davie2007uniqueness}
		Alexander~M Davie.
		\newblock Uniqueness of solutions of stochastic differential equations.
		\newblock {\em International Mathematics Research Notices}, 2007, 2007.
		
		\bibitem[DD09]{dong2009navier}
		Hongjie Dong and Dapeng Du.
		\newblock The {N}avier-{S}tokes equations in the critical {L}ebesgue space.
		\newblock {\em Communications in Mathematical Physics}, 292(3):811--827, 2009.
		
		\bibitem[DPMN92]{da1992compact}
		Giuseppe Da~Prato, Paul Malliavin, and David Nualart.
		\newblock Compact families of {W}iener functionals.
		\newblock {\em {Comptes rendus de l'Acad{\'e}mie des sciences. S{\'e}rie 1,
				Math{\'e}matique}}, 315(12):1287--1291, 1992.
		
		\bibitem[DZ01]{duoandikoetxea2001fourier}
		Javier Duoandikoetxea and Javier~Duoandikoetxea Zuazo.
		\newblock {\em Fourier analysis}, volume~29.
		\newblock American Mathematical Soc., 2001.
		
		\bibitem[ES{\v{S}}03a]{escauriaza2003backward}
		Luis Escauriaza, Gregory Seregin, and Vladimir {\v{S}}ver{\'a}k.
		\newblock Backward uniqueness for parabolic equations.
		\newblock {\em Archive for Rational Mechanics and Analysis}, 169(2):147--157,
		2003.
		
		\bibitem[ES{\v{S}}03b]{escauriaza2003solutions}
		Luis Escauriaza, Gregory Seregin, and Vladimir {\v{S}}ver{\'a}k.
		\newblock ${L}_{3,\infty}$-solutions of the navier-stokes equations and
		backward uniqueness.
		\newblock {\em Russian Mathematical Surveys}, 58(2):211--250, 2003.
		
		\bibitem[Eva10]{evans2010partial}
		Lawrence~C Evans.
		\newblock {\em Partial Differential Equations}.
		\newblock The American Mathematical Society, 2010.
		
		\bibitem[FF11]{fedrizzi2011pathwise}
		Ennio Fedrizzi and Franco Flandoli.
		\newblock Pathwise uniqueness and continuous dependence for {{SDE}}s with
		non-regular drift.
		\newblock {\em {S}tochastics: An International Journal of Probability and
			Stochastic Processes}, 83(03):241--257, 2011.
		
		\bibitem[FF13]{fedrizzi2013noise}
		Ennio Fedrizzi and Franco Flandoli.
		\newblock Noise prevents singularities in linear transport equations.
		\newblock {\em Journal of Functional Analysis}, 264(6):1329--1354, 2013.
		
		\bibitem[FGP10]{flandoli2010well}
		Franco Flandoli, Massimiliano Gubinelli, and Enrico Priola.
		\newblock Well-posedness of the transport equation by stochastic perturbation.
		\newblock {\em Inventiones mathematicae}, 180(1):1--53, 2010.
		
		\bibitem[FJR72]{fabes1972initial}
		Eugene~Barry Fabes, B~Frank Jones, and Nestor~M Riviere.
		\newblock The initial value problem for the {Navier-Stokes} equations with data
		in ${L}^p$.
		\newblock {\em Archive for Rational Mechanics and Analysis}, 45(3):222--240,
		1972.
		
		\bibitem[Gal24]{galeati2025almost}
		Lucio Galeati.
		\newblock Almost-everywhere uniqueness of {L}agrangian trajectories for 3{D}
		{N}avier-{S}tokes revisited.
		\newblock {\em arXiv preprint arXiv:2406.12788}, 2024.
		
		\bibitem[Gig86]{giga1986solutions}
		Yoshikazu Giga.
		\newblock Solutions for semilinear parabolic equations in ${L}^p$ and
		regularity of weak solutions of the {Navier-Stokes} system.
		\newblock {\em Journal of Differential Equations}, 62(2):186--212, 1986.
		
		\bibitem[GP24]{grafner2024weak}
		Lukas Gräfner and Nicolas Perkowski.
		\newblock Weak well-posedness of energy solutions to singular sdes with
		supercritical distributional drift, 2024.
		
		\bibitem[HZ24]{hao2024sdes}
		Zimo Hao and Xicheng Zhang.
		\newblock Sdes with supercritical distributional drifts, 2024.
		
		\bibitem[KM25]{kinzebulatov2025strong}
		Damir Kinzebulatov and Kodjo~Raphael Madou.
		\newblock Strong solutions of {SDE}s with singular (form-bounded) drift via
		{R}oeckner-{Z}hao approach.
		\newblock {\em Stochastics and Dynamics}, 2025+.
		
		\bibitem[KR05]{krylov2005strong}
		{N}icolai~{V} {K}rylov and {M}ichael {R}\"ockner.
		\newblock {S}trong solutions of stochastic equations with singular time
		dependent drift.
		\newblock {\em {P}robability Theory and Related Fields}, 131(2):154--196, 2005.
		
		\bibitem[Kry01]{krylov2001heat}
		Nicolai~V Krylov.
		\newblock The heat equation in ${L}_q ((0, {T}), {L}_p)$-spaces with weights.
		\newblock {\em SIAM Journal on Mathematical Analysis}, 32(5):1117--1141, 2001.
		
		\bibitem[Kry21a]{krylov2021stochastic2}
		N.~V. Krylov.
		\newblock On stochastic equations with drift in ${L}_d$.
		\newblock {\em The Annals of Probability}, 49(5):2371--2398, 2021.
		
		\bibitem[Kry21b]{krylov2021stochastic1}
		N.~V. Krylov.
		\newblock On stochastic {I}t\^o processes with drift in ${L}_d$.
		\newblock {\em Stochastic Processes and their Applications}, 138:1--25, 2021.
		
		\bibitem[Kry21c]{krylov2021strong}
		N.~V. Krylov.
		\newblock On strong solutions of {I}t{\^o}’s equations with $\sigma\in
		{W}^1_d$ and $b\in {L}_d$.
		\newblock {\em The Annals of Probability}, 49(6):3142--3167, 2021.
		
		\bibitem[Kry25]{krylov2025strong}
		N.~V. Krylov.
		\newblock On weak and strong solutions of time inhomogeneous {I}t\^o's
		equations with {VMO} diffusion and {M}orrey drift.
		\newblock {\em Stochastic Process. Appl.}, 179:Paper No. 104505, 23, 2025.
		
		\bibitem[Lad67]{ladyzhenskaya1967uniqueness}
		Olga~Aleksandrovna Ladyzhenskaya.
		\newblock On the uniqueness and on the smoothness of weak solutions of the
		{Navier-Stokes} equations.
		\newblock {\em Zapiski Nauchnykh Seminarov POMI}, 5:169--185, 1967.
		
		\bibitem[Ler34]{leray1934mouvement}
		Jean Leray.
		\newblock Sur le mouvement d'un liquide visqueux emplissant l'espace.
		\newblock {\em Acta Mathematica}, 63:193--248, 1934.
		
		\bibitem[LT21]{lee2017existence}
		Haesung Lee and Gerald Trutnau.
		\newblock Existence, uniqueness and ergodic properties for time-homogeneous
		{I}t\^{o}-{SDE}s with locally integrable drifts and {S}obolev diffusion
		coefficients.
		\newblock {\em Tohoku Math. J. (2)}, 73(2):159--198, 2021.
		
		\bibitem[MBP10]{meyer2010construction}
		Thilo Meyer-Brandis and Frank Proske.
		\newblock Construction of strong solutions of {{SDE}}'s via {M}alliavin
		calculus.
		\newblock {\em Journal of Functional Analysis}, 258(11):3922--3953, 2010.
		
		\bibitem[MNP15]{mohammed2015sobolev}
		Salah-Eldin~A Mohammed, Torstein~K Nilssen, and Frank~N Proske.
		\newblock {S}obolev differentiable stochastic flows for {{SDE}}s with singular
		coefficients: Applications to the transport equation.
		\newblock {\em The Annals of Probability}, 43(3):1535--1576, 2015.
		
		\bibitem[MPMBN{\etalchar{+}}13]{menoukeu2013variational}
		Olivier Menoukeu-Pamen, Thilo Meyer-Brandis, Torstein Nilssen, Frank Proske,
		and Tusheng Zhang.
		\newblock A variational approach to the construction and {M}alliavin
		differentiability of strong solutions of {{SDE}}’s.
		\newblock {\em Mathematische Annalen}, 357(2):761--799, 2013.
		
		\bibitem[Nam20]{nam2020stochastic}
		Kyeongsik Nam.
		\newblock Stochastic differential equations with critical drifts.
		\newblock {\em Stochastic Processes and their Applications}, 130(9):5366--5393,
		2020.
		
		\bibitem[NO15]{neves2015wellposedness}
		Wladimir Neves and Christian Olivera.
		\newblock Well-posedness for stochastic continuity equations with
		{Ladyzhenskaya--Prodi--Serrin} condition.
		\newblock {\em Nonlinear Differential Equations and Applications NoDEA},
		22(5):1247--1258, 2015.
		
		\bibitem[Pro59]{prodi1959teorema}
		Giovanni Prodi.
		\newblock Un teorema di unicita per le equazioni di {Navier-Stokes}.
		\newblock {\em Annali di Matematica Pura ed Applicata}, 48(1):173--182, 1959.
		
		\bibitem[Rez14]{rezakhanlou2014regular}
		Fraydoun Rezakhanlou.
		\newblock Regular flows for diffusions with rough drifts.
		\newblock {\em arXiv preprint arXiv:1405.5856}, 2014.
		
		\bibitem[Rez16]{rezakhanlou2016stochastically}
		Fraydoun Rezakhanlou.
		\newblock Stochastically symplectic maps and their applications to the
		{Navier-Stokes} equation.
		\newblock {\em {Annales de l'Institut Henri Poincare (C) Non Linear Analysis}},
		33(1):1--22, 2016.
		
		\bibitem[RZ23]{rockner2022weak}
		Michael R\"{o}ckner and Guohuan Zhao.
		\newblock S{DE}s with critical time dependent drifts: weak solutions.
		\newblock {\em Bernoulli}, 29(1):757--784, 2023.
		
		\bibitem[Ser62]{serrin1962interior}
		James Serrin.
		\newblock On the interior regularity of weak solutions of the {Navier-Stokes}
		equations.
		\newblock {\em Archive for Rational Mechanics and Analysis}, 9:187--195, 1962.
		
		\bibitem[Sha16]{shaposhnikov2016some}
		AV~Shaposhnikov.
		\newblock Some remarks on {D}avie's uniqueness theorem.
		\newblock {\em Proceedings of the Edinburgh Mathematical Society},
		59(4):1019--1035, 2016.
		
		\bibitem[Sob77]{sobolevskii1977fractional}
		Pavel~Evseyevich Sobolevskii.
		\newblock Fractional powers of coercive-positive sums of operators.
		\newblock {\em Siberian Mathematical Journal}, 18(3):454--469, 1977.
		
		\bibitem[Ver80]{veretennikov1980strong2}
		Alexander~Yur'evich Veretennikov.
		\newblock On strong solutions and explicit formulas for solutions of stochastic
		integral equations.
		\newblock {\em Matematicheskii Sbornik}, 153(3):434--452, 1980.
		
		\bibitem[VK76]{veretennikov1976explicit}
		A~Ju Veretennikov and Nicolai~V Krylov.
		\newblock On explicit formulas for solutions of stochastic equations.
		\newblock {\em Mathematics of the USSR-Sbornik}, 29(2):239--256, 1976.
		
		\bibitem[XXZZ20]{xia2020lqlp}
		Pengcheng Xia, Longjie Xie, Xicheng Zhang, and Guohuan Zhao.
		\newblock ${L}^q({L}^p)$-theory of stochastic differential equations.
		\newblock {\em Stochastic Processes and their Applications}, 130(8):5188--5211,
		2020.
		
		\bibitem[Zha05]{zhang2005strong}
		Xicheng Zhang.
		\newblock Strong solutions of {{SDE}}s with singular drift and {S}obolev
		diffusion coefficients.
		\newblock {\em Stochastic Processes and their Applications},
		115(11):1805--1818, 2005.
		
		\bibitem[{Z}ha10]{zhang2010astochastic}
		{X}icheng {Z}hang.
		\newblock {A} stochastic representation for backward incompressible
		{{N}}avier-{S}tokes equations.
		\newblock {\em {P}robability Theory and Related Fields}, 148(1-2):305--332,
		2010.
		
		\bibitem[{Z}ha11]{zhang2011stochastic}
		{X}icheng {Z}hang.
		\newblock {S}tochastic homeomorphism flows of {S}{D}{E}s with singular drifts
		and {S}obolev diffusion coefficients.
		\newblock {\em {E}lectronic {J}ournal of {P}robability}, 16:1096--1116, 2011.
		
		\bibitem[{Z}ha16]{zhang2016stochastic}
		{X}icheng {Z}hang.
		\newblock {S}tochastic differential equations with {S}obolev diffusion and
		singular drift and applications.
		\newblock {\em {T}he {A}nnals of {A}pplied {P}robability}, 26(5):2697--2732,
		2016.
		
		\bibitem[Zvo74]{zvonkin1974transformation}
		Alexander~K Zvonkin.
		\newblock A transformation of the phase space of a diffusion process that
		removes the drift.
		\newblock {\em Mathematics of the USSR-Sbornik}, 22(1):129, 1974.
		
		\bibitem[ZZ18]{zhang2018singular}
		Xicheng Zhang and Guohuan Zhao.
		\newblock {S}ingular {B}rownian diffusion processes.
		\newblock {\em {C}ommunications in {M}athematics and {S}tatistics},
		6(4):533--581, 2018.
		
		\bibitem[ZZ21]{zhang2020stochastic}
		Xicheng Zhang and Guohuan Zhao.
		\newblock Stochastic {L}agrangian path for {L}eray’s solutions of 3{D}
		{N}avier--{S}tokes equations.
		\newblock {\em Communications in Mathematical Physics}, 381(2):491--525, 2021.
		
	\end{thebibliography}
\end{document}